\documentclass[reqno,11pt,a4paper]{amsart}

\usepackage[utf8]{inputenc}
\usepackage[T1]{fontenc}
\usepackage{ae}
\usepackage{units}
\usepackage{icomma}
\usepackage{color}
\usepackage{graphicx}
\usepackage{cancel}
\usepackage{fullpage}
\usepackage{upgreek}
\usepackage{algorithm}
\usepackage{algorithmic}
\usepackage{type1cm}
\usepackage{eso-pic}
\usepackage[breaklinks,pdfpagelabels=false]{hyperref}
\usepackage{verbatim}
\usepackage{tikz}
\usetikzlibrary{arrows}

\usepackage{caption}
\renewcommand{\epsilon}{\varepsilon}

\newtheorem{theorem}{Theorem}[section]
\newtheorem{lemma}[theorem]{Lemma}
\newtheorem{proposition}[theorem]{Proposition}
\newtheorem{corollary}[theorem]{Corollary}

\newtheorem{conjecture}[theorem]{Conjecture}
\newtheorem{remark}[theorem]{Remark}
\theoremstyle{definition}
\newtheorem{definition}[theorem]{Definition}

\renewcommand{\qed}{\hfill \mbox{\raggedright \rule{0.1in}{0.1in}}}

\numberwithin{equation}{section}
\numberwithin{theorem}{section}

\makeatletter

\def\@settitle{\begin{center}\baselineskip14\p@\relax\bfseries\Large\@title\end{center}}

\makeatother

\begin{document}

\title[]
{On $k$-sum-free subsets of $\{-n,\,\dots,\,n\}$}



      \date{\today}


      \author{Peter Hegarty and Damiano Ognissanti} \address{Department of Mathematical Sciences, Chalmers University of Technology and University of Gothenburg, 41296 Gothenburg, Sweden} 
\email{hegarty@chalmers.se, damiano@chalmers.se}







\begin{abstract}
  For each $k \geq 1$, we determine the maximum asymptotic density, as
  $n \rightarrow \infty$, of a $k$-sum-free subset of $\{-n,\,\dots,\,n\}$,
  that is, of a set containing no non-trivial solutions to the equation
  $kx=y+z$. These are analogues of known results for subsets of
  $\{1,\,\dots,\,n\}$, and are genuinely new for $k \geq 4$: for all such $k$,
  the maximum density is strictly less than for subsets of
  $\{1,\,\dots,\,n\}$, so not achieved by $\pm S$,
  where $S$ is a maximum $k$-sum-free subset of $\{1,\,\dots,\,n\}$. Further,
  we obtain the approximate structure of all extremal $k$-sum-free sets in
  $\{-n,\,\dots,\,n\}$, up to $O(1)$ elements for $k \geq 5$, and up to $o(n)$
  elements for $k=4$, where the extremal sets have a different structure. We
  also present continuous analogues of our results, for measurable $k$-sum-free
  subsets of the real interval $[-1,\,1]$. Here, $k=3$ remains open and
  numerical evidence suggests a different structure for the extremal sets than
  when $k \geq 4$. Section \ref{sect:dilates} contains a result on sums of
  dilates which may be of independent interest. 
  \end{abstract}


\maketitle

\setcounter{equation}{0}

\setcounter{equation}{0}
\setcounter{section}{0}


\section{{\bf Introduction}}\label{sect:intro}

Let $a,\,b,\,c$ be positive integers, $G$ an abelian group. A solution in $G$
to the linear equation $\mathcal{L} = \mathcal{L}_{a,\,b,\,c}: ax=by+cz$ is
called trivial if
$x=y=z$ and non-trivial otherwise. A subset $S$ of $G\backslash \{0\}$
is said to be \emph{$(a,\,b,\,c)$-free} if there is no non-trivial solution in
$S$ to $\mathcal{L}$. If $G$ is
torsion-free, then trivial solutions for $x \neq 0$ can only
occur when $a=b+c$, in which case $\mathcal{L}$ is said to be
\emph{(translation) invariant}. Otherwise, $\mathcal{L}$ is said to be
\emph{non-invariant}. 

The classical setting in which $(a,\,b,\,c)$-free sets have been studied is
where $G = \mathbb{Z}$ and $S \subseteq [1,\,n]$ is a finite set of positive
integers. The basic question of interest is how dense $S$ can be, which means
investigating
\begin{eqnarray}\label{eq:density}
  D_n = D_n (a,\,b,\,c) =  \frac{1}{n} \times \max \{|S| : S \subseteq [1,\,n] \; {\hbox{and $S$ is $(a,\,b,\,c)$-free}} \}, \nonumber \\
  D = D(a,\,b,\,c) = \limsup_{n \rightarrow \infty} D_n.
    \end{eqnarray}
    The basic dichotomy is between invariant and non-invariant equations. If
    the equation is invariant, then $D = 0$. This is highly non-trivial. The
    special case $(a,\,b,\,c) = (2,\,1,\,1)$ is Roth's theorem and his method
    can be applied more generally: see \cite{R1}. 
    On the other hand, it is easy to see that $D > 0$ for any
    non-invariant equation, since for each $n$ we can take either
    \begin{itemize}
    \item
      $S = \{s \in [1,\,n] : s \equiv 1 \; ({\hbox{mod $t$}}) \}$, for any $t$
      that doesn't divide $l = |a-(b+c)|$, or
    \item
      $S = (\varepsilon n, \, n]$, where
      $\varepsilon = \frac{m}{M}$ and $m, \, M$ are, respectively, the min and
      max of $\{a,\,b+c\}$. Throughout this paper, except for Section
      \ref{sect:cont}, interval notation
      denotes sets of integers: thus $(\alpha,\,\beta) = \{x \in \mathbb{Z}: \alpha < x < \beta\}$, and similarly for (half)-closed intervals. 
    \end{itemize}

    This paper is concerned with non-invariant equations. Since we know $D > 0$
    always, the next obvious problem is to compute the function
    $D:\mathbb{N}^3 \backslash \mathbb{H} \rightarrow (0,\,1]$, where
    $\mathbb{H}$ is the hyperplane $a=b+c$. This problem is wide open in general.
    An indication of its difficulty is given by the fact that it is not even
    known in general if the limsup in (\ref{eq:density}) is actually a limit{\footnote{We are not completely sure that this remains an open problem. It was originally posed in \cite{R2}. At the end of the paper \cite{CS} from 2011, the authors say they will present a solution of this problem in a forthcoming paper. We have not been able to locate any such solution in the post-2011 literature, however.}}.
    There are, however, some cases where we know the answer and the most famous
    of these is probably $(a,\,b,\,c) = (k,\,1,\,1)$, for some
    $k \in \mathbb{N}$, so-called \emph{$k$-sum-free} sets. Let
    $D^k = D(k,\,1,\,1)$. Then the following are known:
    \begin{itemize}
    \item
      $k=1$: $D^1 = \frac{1}{2}$. For sufficiently large $n$, there are either
      two or three maximum $1$-free subsets of $[1,\,n]$, the odd numbers and
      \begin{itemize}
        \item
          the interval $\left( \frac{n}{2}, \, n\right]$, when $n$ is odd.
        \item
          both $\left( \frac{n}{2}, \, n\right]$ and
          $\left[ \frac{n}{2}, \, n\right)$, when $n$ is even.
        \end{itemize}
        These facts follow immediately from the results in \cite{D}.
    \item
      $k=2$: The equation is invariant, $D^2 = 0$.
    \item
      $k=3$: $D^{3} = \frac{1}{2}$. The odd numbers form the unique extremal
      subset of $[1,\,n]$ for all sufficiently large $n$.
      This was proven in \cite{CG}.
    \item
      $k \geq 4$. For all such $k$ one has
      \begin{eqnarray*}
        D^k = \frac{(k-2)(k+s)}{k^2 - 2}, \;\;\; {\hbox{where}} \; s = \frac{8}{k(k^4 - 2k^2 - 4)}.
    \end{eqnarray*}
    Maximum asymptotic density is achieved by the sets
    \begin{eqnarray}\label{eq:Splus}
      S^{+}_{k,\,n} = \left( \frac{2}{k} \varepsilon_2 n, \, \varepsilon_2 n \right] \cup \left( \frac{2}{k} \varepsilon_1 n, \, \varepsilon_1 n \right] \cup \left( \frac{2}{k} n, \, n \right], \;\;\; {\hbox{where}} \; \varepsilon_2 = \frac{ks}{2}, \; \varepsilon_1 = \frac{2}{k^2}(1 + \varepsilon_2).
    \end{eqnarray}
    Moreover, these are essentially the only extremal sets: there are numbers
    $C(k)$ such that, if $S$ is any
    maximum $k$-sum-free subset of $[1,\,n]$, then
    $|S \Delta S^{+}_{k,\,n}| \leq C(k)$. All of this was proven in \cite{H1}.
  \end{itemize}

  In the context of finite $\mathcal{L}$-free sets, it may be considered more
  natural to consider subsets of $\mathbb{Z}_n$, the integers modulo $n$. To
  avoid technicalities caused by the presence of torsion, attention is
  normally confined to when $n=p$, a prime. We can make
  definitions analogous to (\ref{eq:density}):
\begin{eqnarray}\label{eq:density2}
  \delta_n = \delta_n (a,\,b,\,c) =  \frac{1}{n} \times \max \{|S| : S \subseteq \mathbb{Z}_n \; {\hbox{and $S$ is $(a,\,b,\,c)$-free}} \}, \nonumber \\
  \delta = \delta(a,\,b,\,c) = \limsup_{p \rightarrow \infty, \, {\hbox{$p$ prime}}} \delta_p.
    \end{eqnarray}
    The basic dichotomy between invariant and non-invariant equations is easily
    seen to remain. In contrast to the positive integer setting, it is known
    that the limsup in (\ref{eq:density2}) is indeed a limit, and that the
    limit in question is the supremum of the measures of Lebesgue measurable
    $(a,\,b,\,c)$-free subsets of $\mathbb{R}/\mathbb{Z}$. This was proven in
    \cite{CS}. The relationship between the discrete and continuous problems is
    more complicated in the integer setting, see Section \ref{sect:cont}
    below.
    
    On the other hand, the actual determination of maximum densities, for
    non-invariant $\mathcal{L}$, seems
    to be more difficult in the modular setting. Neither of the basic ideas
    mentioned above for
    constructing positive density $\mathcal{L}$-free sets
    transfers straightforwardly from integers in $[1,\,n]$ to integers
    modulo $p$. Indeed, a construction based on congruence classes doesn't even
    make sense in the new context. A construction based on ``edge intervals''
    won't work as well as before, because there are no ``edges'' in the modular
    setting. 
    While it's trivial that $\delta_p \leq D_p$ for
    any equation $\mathcal{L}$ and any prime $p$, these considerations indicate
    why the problem of computing $\delta_p, \, \delta$ in the latter setting
    might in general be more difficult. For the case of $k$-sum-free sets,
    this indeed seems to be the case. The following is what is currently known:
  \begin{itemize}
  \item
    $k=1$: The Cauchy-Davenport theorem immediately implies that, for any
    fixed $k \neq 2$ and prime $p > k$, a $k$-sum-free
    subset of $\mathbb{Z}_p$ can have size at most
    $\lfloor \frac{p+1}{3}\rfloor$. Equality is achieved for $k=1$ and
    every $p$ by the
    middle-third $\left( \frac{p}{3}, \, \frac{2p}{3} \right]$. Hence, the
    problem is no more difficult than for natural numbers, but the maximum
    asymptotic density drops from $\frac{1}{2}$ to $\frac{1}{3}$.
    \item
      $k=2$: The equation is invariant, and the asymptotic density can't be
      non-zero. This is just as for sets of natural numbers.
    \item
      $k \geq 3$: Here is where things seem to get much more difficult. The
      maximum asymptotic density is not known for any $k$, but is known to be
      strictly less than $\frac{1}{3}$ for every $k$. For the state of the art
      on this problem, see \cite{C}.
    \end{itemize}

    Our idea for this paper was to look at $k$-sum-free subsets of
    $[-n,\,n] \cap \mathbb{Z}$. The definitions analogous to (\ref{eq:density})
    and (\ref{eq:density2}) are
\begin{eqnarray}\label{eq:density3}
  \partial_n = \partial_n (a,\,b,\,c) =  \frac{1}{2n+1} \times \max \{|S| : S \subseteq [-n,\,n] \; {\hbox{and $S$ is $(a,\,b,\,c)$-free}} \}, \nonumber \\ \partial = \partial(a,\,b,\,c) = \limsup_{n \rightarrow \infty} \partial_n.
    \end{eqnarray}
    As far as we're aware, there is no previous literature on extremal
    $(a,\,b,\,c)$-free sets in this setting. It is somewhere
    ``in between'' the two traditional ones of
    $[1,\,n] \cap \mathbb{Z}$ and $\mathbb{Z}_n$ - indeed, it is trivial that,
    for any equation $\mathcal{L}$ and $n > 0$,
    \begin{eqnarray}\label{eq:trivineq}
      \delta_{2n+1} \leq \partial_n < D_n,
    \end{eqnarray}
    where the strict inequality between $\partial_n$ and $D_n$ is because an
    $\mathcal{L}$-free set can't contain zero. In particular, asymptotic
    density must still be zero for invariant $\mathcal{L}$, so we can still
    focus on
    non-invariant equations.

    Since we're still working in
    the group $\mathbb{Z}$, we 
    expect behaviour similar to that
    for positive integers. For example, if a maximum $\mathcal{L}$-free
    $S \subseteq [1,\,n]$ consists of a symmetric set of congruence classes,
    that is $S \equiv -S ({\hbox{mod $m$}})$ for some $m$, then $\pm S$
    will automatically be maximum $\mathcal{L}$-free in $[-n,\,n]$.
    The symmetric interval also still has ``edges'' which can be exploited to
    construct other types of dense
    $\mathcal{L}$-free sets. On the other
    hand, there is now the possibility of cancellation between positive and
    negative integers, which mimics a feature of the $\mathbb{Z}_n$ setting (as
    in any abelian group). Intuitively, then, one should in general expect
    extremal $\mathcal{L}$-free sets in this setting to look like those in
    the positive integer setting, though perhaps with strictly smaller
    achievable asymptotic density when the extremal sets contain
    ``edge intervals''.
    \\
    \par We will confirm this intuition for $k$-sum-free sets. Let
    $\partial^k = \partial(k,\,1,\,1)$. We will establish the following picture:
    \begin{itemize}
    \item
      $k=1$: $\partial^1 = \frac{1}{2}$. This is a trivial extension of the
      result for positive integers. We can still take either the odd numbers in
      $[-n,\,n]$, or the union of two edge intervals
      $\pm \left( \frac{n}{2}, \, n \right)$ plus some combination of
      endpoints. Moreover, these are the
      only extremal sets for large enough $n$.
    \item
      $k=2$: The equation is invariant, so $\partial^2 = 0$ just as in the
      other two settings.
    \item
      $k=3$: $\partial^3 = \frac{1}{2}$. It is a trivial consequence of the
      result for positive integers that the odd numbers in $[-n,\,n]$ form the
      unique extremal set for all $n \gg 0$.
    \item
      $k \geq 4$: This is where things get interesting. We will show that
      $\partial^k < D^k$ for all $k \geq 4$. Moreover, the extremal sets look
      different for $k \geq 5$ than for $k=4$. For $k \geq 5$ we cannot do
      better than $\pm I$ where $I$ is an optimal edge interval. But for
      $k = 4$, a set consisting of three intervals, two edge intervals and a
      suitably chosen third interval entirely on one side of zero, achieves
      optimal asymptotic density. 
    \end{itemize}

    Let us now formulate our main results more precisely. For positive integers
    $k, \, n$ with $k \geq 4$, let
    \begin{eqnarray}\label{eq:S}
      S_{k,\,n} = \pm \left( \frac{2n}{k}, \, n \right].
    \end{eqnarray}
    Let also
    \begin{eqnarray}\label{eq:S4}
      \mathcal{S}_{4,\,n} = \pm \left( \frac{n}{2}, \, n \right] \cup \left[ \frac{n}{8}, \, \frac{5n}{32} \right] = S_{4,\,n} \cup \left[ \frac{n}{8}, \, \frac{5n}{32} \right].
    \end{eqnarray}
    It is easy to check that $S_{k,\,n}$ is $k$-sum-free and that
    $\mathcal{S}_{4,\,n}$ is $4$-sum-free. Note that, as $n \rightarrow \infty$,
    \begin{eqnarray}\label{eq:Slim}
      \frac{|S_{k,\,n}|}{2n+1} \rightarrow \frac{k-2}{k}, \;\;\;\; \frac{|\mathcal{S}_{4,\,n}|}{2n+1} \rightarrow \frac{33}{64}.
    \end{eqnarray}
    Our two main results are as follows. Recall that $A \Delta B$ denotes the
    symmetric difference of sets $A$ and $B$, and that $-A = \{-a: a \in A\}$. 

    \begin{theorem}\label{thm:A}
      There is an absolute constant $C$ such that, for all $k \geq 5$ and
      $n \geq 1$, if $S$ is a $k$-sum-free subset of
      $[-n,\,n]$ of maximum size, then $|S \Delta S_{k,\,n}| \leq C$. In
      particular, $\partial^k = \frac{k-2}{k}$. 
    \end{theorem}

    \begin{theorem}\label{thm:B}
      There is a function $f(n)$ such that $f(n) = o(n)$ and, if $S$ is a
      $4$-sum-free subset of $[-n,\,n]$ of maximum size, then
      $\min \{|S \Delta \mathcal{S}_{4,\,n}|, \, |S \Delta -\mathcal{S}_{4,\,n}| \} \leq f(n)$. In particular, $\partial^4 = \frac{33}{64}$.
    \end{theorem}

    Computations suggest, but we were unable to prove, the
    following strengthening of Theorem \ref{thm:B}:

    \begin{conjecture}\label{conj:C}
      There is a constant $C$, independent of $n$, such that if $S$ is a
      $4$-sum-free subset of $[-n,\,n]$ of maximum size, then
      $\min \{|S \Delta \mathcal{S}_{4,\,n}|, \, |S \Delta -\mathcal{S}_{4,\,n}| \} \leq C$.
    \end{conjecture}

    Our computations have not yielded any example where $C > 6$.
    \\
    \par The rest of the paper is organized as follows. Section
    \ref{sect:lemmas} introduces some notation and terminology, and the basic
    strategy behind the proofs of our theorems.
    The overall idea here is that, in attempting to
    construct a $k$-sum-free subset $S$ of $[-n,\,n]$ which is larger than the
    set $S_{k,\,n}$, or $\mathcal{S}_{4,\,n}$ for $k=4$, we will inevitably have
    to ``trade off'' elements from other parts of the interval $[-n,\,n]$ with
    those from the edge intervals $\pm \left( \frac{2n}{k}, \, n \right]$. This
    trade-off is always bad, but how bad depends on where the other elements in
    $S$ are taken from. The problem is that the parts of the edge intervals
    that are being ``traded off'' against other parts of $[-n,\,n]$ are not a
    priori disjoint. Section \ref{sect:lemmas} contains just enough material
    necessary to prove Theorem \ref{thm:A} for $k \geq 6$, and
    this is completed easily in Section \ref{sect:ksix}. For
    $k = 5$ we need to be somewhat more careful, but the proof is still
    relatively easy: we complete the proof of Theorem \ref{thm:A} in Section
    \ref{sect:kfive}. Sections \ref{sect:dilates}-\ref{sect:stability} are
    devoted to the proof of Theorem \ref{thm:B}. For $k=4$, a considerably more
    complicated leveraging of the basic ideas from previous sections is
    necessary: there turn out to be many ways in which elements from different
    parts of $[-n,\,n]$ can be ``traded off'' and we eventually find
    suitable ``combinations'' of these trade-offs which allow us to prove the
    theorem. In Section \ref{sect:dilates}, we prove a general result on the
    sizes of sets of the form $4A - A$, which may be of independent interest.
    It is a kind of ``one-sided'' version of (a special case of) a well-known
    result of Shakan
    \cite{S} on sums of dilates. In Section \ref{sect:kfour}, we prove a
    weaker version of Theorem \ref{thm:B}, namely 
    that $|S| \leq \frac{33}{32} n + o(n)$, for any $4$-sum-free set
    $S \subseteq [-n,\,n]$. This is by far the longest section with the most
    intricate analysis of trade-offs - indeed, some of the most technical are
    relegated to an Appendix.
    In Section \ref{sect:stability} we go through the
    further work that needs to be done to prove Theorem \ref{thm:B} in full.
    The arguments from Section \ref{sect:kfour} rule out all but essentially
    one type of configuration as candidates for maximum $4$-free subsets of
    $[-n,\,n]$, and in Section \ref{sect:stability} we prove that such a
    configuration, if maximum, must be very close to $\mathcal{S}_{4,\,n}$. 
    Section \ref{sect:cont} discusses natural continuous
    analogues of the main results, where we consider $k$-sum-free measurable
    subsets of the real interval $[-1,\,1]$. Because the extremal sets in the
    discrete setting are, up to $o(1)$ in density, unions of finitely many
    intervals, a restricted ``transference'' principle, different from and
    simpler
    than that employed in \cite{CS} for the modular setting, establishes a
    similar result for measurable subsets of $[-1,\,1]$.
    However, for $k=3$ the symmetric continuous problem
    remains open - it was resolved in \cite{PR} for subsets of $[0,\,1]$.
    We present a conjecture which is supported by numerical evidence, but where
    indications so far are that it might require arguments at least as
    intricate as those in earlier sections of this paper.
    Section \ref{sect:AI} is our AI statement, describing
    the role of AI in producing the results of the paper and the contributions
    of the human authors.
    
    \vspace{3mm}

\section{{\bf Basic Tools}}\label{sect:lemmas}

From now until the end of Section \ref{sect:stability}, the following notation
applies. The letters $k$ and $n$ always denote positive integers, such that
$k \geq 4$. A reminder that we use interval notation for sets of integers: thus
$(\alpha,\,\beta)$ denotes $\{x \in \mathbb{Z}: \alpha < x < \beta\}$, and
similarly for (half)-closed
intervals. Note that, for any $\alpha < \beta$, there are at most
$(\beta-\alpha) + 1$ integers
in such an interval. The letter $S$ always denotes a $k$-sum-free subset of
$[-n,\,n]$. When the value of $k$ is understood, we just say $S$ is free.
Recall that $0 \not\in S$. We denote $P = P_S = S \cap [1,\,n]$ and
$N = N_S = -(S \cap [-n,\, -1])$. Thus $S = P \sqcup -N$. We will repeatedly
make use of the following elementary properties of $k$-sum-free sets:

\begin{lemma}\label{lem:elem}
  Suppose $S$ is $k$-sum-free. Then
  \par (i) Both $P$ and $N$ are $k$-sum-free.
  \par (ii) Let $X \in \{P,\,N\}$. If $a,\,b \in X$, then $ka-b \not\in X$.
  \par (iii) If $p \in P$ and $n \in N$, then $kp + n \not\in P$.
  \par (iv) If $n \in N$ and $p \in P$, then $kn + p \not\in N$.
  \par (v) If $x \in S$, then $kx \neq p-n$ for any $p \in P$ and $n \in N$.
\end{lemma}

\begin{proof}
  All these statements are immediate consequences of the definition of
  $k$-sum-freeness.
\end{proof}

Set $l_k = \min \{k, \, 6\}$. We define subsets of $[1,\,n]$ as follows:
\begin{eqnarray}\label{eq:LMQT}
  L = \left[1, \, \frac{l_k-2}{k^2} n \right), \;\; M = \left[ \frac{l_k-2}{k^2} n, \, \frac{4}{k^2} n \right], \;\; Q = \left( \frac{4}{k^2} n, \, \frac{2}{k} n \right], \;\; T = \left( \frac{2}{k} n, \, n \right].
\end{eqnarray}
Thus, the intervals $L$ (``lower''), $M$ (``middle''), $Q$ and $T$ (``top'')
partition $[1,\,n]$. But $M$ consists of at
most a single point for $k \geq 6$, so it really only plays a role for
$k \in \{4,\,5\}$.

Note that $S_{k,\,n} = \pm T$. The basic strategy behind our proofs is that
presence in a $k$-sum-free $S$ of elements from any of the intervals
$\pm L, \, \pm M$ or $\pm Q$ must be traded off against the presence of
elements of $\pm T$. The trade-off is bad in each case, but how bad depends on
the specific interval. To complete a proof, we need to find enough pairwise
disjoint parts of $\pm T$ that must be traded off against (in absolute value)
smaller numbers. Moreover, for $k = 4$, we have the problem that a part of
$M = \left[ \frac{n}{8}, \, \frac{n}{4} \right]$, namely
$\left[\frac{n}{8}, \, \frac{5n}{32} \right]$, can be added to $\pm T$ without
any trade-off.
\\
\par In this section, we prove just two lemmas about ``trade-offs'', which will
be used in the subsequent analysis for all values of $k$. Part (ii) of the
first of these can be thought of as a price paid for the presence of elements
from $\pm L$, rather than a trade-off between $L$ and $T$ necessarily. The
second lemma asserts an explicit trade-off between $Q$ and $T$. It will turn
out that these already suffice to prove Theorem \ref{thm:A} for $k \geq 6$. In
Sections \ref{sect:kfive} and \ref{sect:kfour}, we will derive many other
trade-offs for the intervals $L, \, M$ and $Q$ versus $T$.

\begin{lemma}\label{lem:refred}
  (i) Let $S \subseteq [-n,\,n]$ be $k$-sum-free and $x \in S$ with
  $|x| \leq \frac{2n}{k}$. If $x > 0$, set
  $I = [kx-n,\,n]$. If $x < 0$, set $I = [-n,\,kx+n]$. Then
  $|S \cap I| \leq n - \frac{k|x|}{2} + \frac{1}{2}$.
  \par (ii) Let $c \geq 0$ be a constant and suppose there is some $x \in S$
  with $|x| \leq \frac{2(k-4)}{k^2} \, n - c$. Then $|S| \leq \left( 2 - \frac{4}{k} \right) n - \frac{kc-1}{2}$.
\end{lemma}

\begin{proof} (i) Note that $|I| = 2n - k|x| + 1$. The map
  $\varepsilon: t \mapsto kx - t$ is an involution of $I$ and $t$ and
  $\varepsilon(t)$ cannot both belong to $S$, for
  any $t$. Thus $|S \cap I| \leq \lfloor \frac{|I|}{2} \rfloor \leq n - \frac{k|x|}{2} + \frac{1}{2}$.
  \par (ii) Without loss of generality, $x > 0$. Then, using part (i),
  \begin{eqnarray*}
    |S| = |S \cap [-n,\,kx-n)| + |S \cap I| \leq kx + \left( n - \frac{kx}{2} + \frac{1}{2} \right) \\ = n + \frac{kx}{2} + \frac{1}{2} \leq n + \frac{k}{2} \left( \frac{2(k-4)}{k^2} n - c \right) + \frac{1}{2} = \left( 2 - \frac{4}{k} \right) n - \frac{kc-1}{2}.
    \end{eqnarray*}
  \end{proof}

  \begin{lemma}\label{lem:Qtrade}
    Let $S \subseteq [-n,\,n]$ be $k$-sum-free. Then
    \begin{eqnarray}\label{eq:Qtrade}
      |Q \cap S| \leq \frac{4}{k} |T \cap S^c| + 2.
    \end{eqnarray}
  \end{lemma}

  \begin{proof}
    Let $t_0 = \min T = \lfloor \frac{2n}{k} \rfloor + 1$ and
    $g = |T \cap S^c|$. For $x \in Q \cap S$, let $t = kx$ and
    $R = R_x = \{y \in T: t-y \in T\}$. Then $y$ and
    $t-y$ cannot both be in $S$ for any $y \in R$. Thus,
    \begin{eqnarray}\label{eq:gR}
      g \geq \frac{|R_x|}{2}, \;\; {\hbox{for any $x \in Q \cap S$}}.
      \end{eqnarray}
      The elements of $R_x$ are the integers $y$ for which
      \begin{eqnarray}\label{eq:Qspan}
        \max \{t_0, \, t-n \} \leq y \leq \min \{n,\,t-t_0\}.
      \end{eqnarray}
      Thus,
      \begin{eqnarray}\label{eq:Rx}
        |R_x| = \left\{ \begin{array}{lr} t - 2t_0 + 1, & {\hbox{if $t \leq n + t_0$}}, \\ 2n-t+1, & {\hbox{if $t \geq n + t_0$}}. \end{array} \right.
      \end{eqnarray}
      It follows from (\ref{eq:gR}) and (\ref{eq:Rx}) that
      \begin{eqnarray*}
        t \leq n + t_0 \; \Rightarrow \; t \leq 2t_0 + 2g - 1, \\
        t \geq n + t_0 \; \Rightarrow \; t \geq 2n - 2g + 1.
      \end{eqnarray*}
      Let $Q_1 = \{x \in Q \cap S : kx \leq 2t_0 + 2g - 1\}$ and $Q_2 = \{x \in Q \cap S: kx \geq 2n-2g + 1\}$. Thus, $|Q \cap S| \leq |Q_1| + |Q_2|$.  

      If $x \in Q_1$, then
      \begin{eqnarray*}
        kx \leq 2t_0 + 2g - 1 \leq \frac{4n}{k} + 2g + 1 \; \Rightarrow x \leq \frac{4n}{k^2} + \frac{2g+1}{k}.
      \end{eqnarray*}
      Hence, $x \in \left( \frac{4n}{k^2}, \, \frac{4n}{k^2} + \frac{2g+1}{k} \right]$ and so $|Q_1| \leq \frac{2g+1}{k} + 1$.

      Similarly, if $x \in Q_2$, then $x \in \left[ \frac{2n}{k} - \frac{2g-1}{k}, \, \frac{2n}{k} \right]$ and so $|Q_2| \leq \frac{2g-1}{k} + 1$.
      Adding, we get
      \begin{eqnarray*}
        |Q \cap S| \leq |Q_1| + |Q_2| \leq \left( \frac{2g+1}{k} + 1 \right) + \left( \frac{2g-1}{k} + 1 \right) = \frac{4g}{k} + 2, \;\; {\hbox{as desired.}}
        \end{eqnarray*}
  \end{proof}
  
\vspace{3mm}

\section{{\bf Proof of Theorem \ref{thm:A} for $k \geq 6$}}\label{sect:ksix}

Note that $S_{k,\,n} = \pm T$ and  
\begin{eqnarray}\label{eq:maxsize}
  \left( 2 - \frac{4}{k} \right) n \leq |S_{k,\,n}| = 2|T| < \left( 2 - \frac{4}{k} \right) n + 2.
\end{eqnarray}
Fix $n$ and $k \geq 6$. Let $S \subseteq [-n,\,n]$ be a $k$-sum-free set for
which $|S| \geq |S_{k,\,n}|$. Consider first just the positive part $P$ of $S$,
for simplicity. Lemma \ref{lem:refred}(ii) and (\ref{eq:maxsize}) imply there
is at most one $x \in P$ such that $x < \frac{2(k-4)}{k^2} n$. Since
$k \geq 6$, this means that $|P \cap L| \leq 1$. Since also $|M| \leq 1$ for
$k \geq 6$, this means $|P \cap (L \cup M)| \leq 2$. Using Lemma
\ref{lem:Qtrade},
\begin{eqnarray*}
  |P| \leq 2 + |P \cap Q| + |P \cap T| \leq 2 + \frac{4}{k} |P^c \cap T| + 2 + |P \cap T| = |T| + 4 - \frac{k-4}{k} |P^c \cap T|.
\end{eqnarray*}
A similar analysis applies to $N$ and so
\begin{eqnarray*}
  |S| \leq 2|T| + 8 - \frac{k-4}{k} (|P^c \cap T| + |N^c \cap T|).
\end{eqnarray*}
Since $|S| \geq 2|T|$, it follows that
$|P^c \cap T| + |N^c \cap T| \leq \frac{8k}{k-4} \leq 24$.
Lemma \ref{lem:Qtrade} then gives
$|P \cap Q| + |N \cap Q| \leq \frac{4}{k} \cdot 24 + 4 \leq 20$. Thus,
$|S \backslash (\pm T)| \leq 2 \cdot 2 + 20 = 24$ and so
$|S \Delta (\pm T)| \leq 48$.
\qed

\vspace{3mm}

\section{{\bf Proof of Theorem \ref{thm:A} for $k = 5$}}\label{sect:kfive}

In remaining sections we will use results on $k$-sum-free subsets of $[1,\,n]$
proven in \cite{H1}. The first of these is

\begin{lemma}\label{lem:H1a}
  Let $A \subseteq [1,\,n]$ be a $k$-sum-free set with smallest element $s$.
  Then
  \begin{eqnarray*}
    |A| \leq |A \cap [1,\,r]| + \lceil \frac{k-2}{k} \, n \rceil, \;\;\; {\hbox{where $r = \lfloor \frac{\lfloor 2n/k \rfloor + s}{k} \rfloor$}}.
  \end{eqnarray*}
\end{lemma}

\begin{proof}
  This is inequality (2) in \cite{H1}.
\end{proof}

Now let $S$ be a $5$-sum-free subset of $[-n,\,n]$ such that
$|S| \geq |S_{5,\,n}|$. Unlike in the previous section, we will no longer try
to compute explicit constants. Firstly, applying Lemma \ref{lem:refred} with
$k=5$ implies that at most $O(1)$ elements $x \in S$ satisfy
$|x| \leq \frac{2n}{25}$. Without loss of generality, we can thus assume $|x| \geq \frac{2n}{25}$ for every $x \in S$ and that instead $|S| \geq |S_{5,\,n}| - O(1)$. Take $\alpha = \frac{21}{125}$ and $\beta = \frac{33}{100}$. There is
some flexibility in the choice of $\alpha$ and $\beta$, we choose specific
values for convenience.
\\
\\
{\sc Case 1:} At least one of $P$ and $N$ meets the interval $[\alpha n, \, \beta n]$.
\\
\par This is the ``easy'' case. Without loss of generality, suppose $P$ meets
this interval. We will show that
\begin{eqnarray}\label{eq:Pab}
  |P| \leq \frac{583}{1000} n + O(1).
\end{eqnarray}
For the moment, suppose this is true. Taking $s = \frac{2n}{25}$ in Lemma
\ref{lem:H1a} we get
\begin{eqnarray}\label{eq:Nstar}
  |N| \leq (r-s+1) + \lceil \frac{3n}{5} \rceil \leq \frac{77}{125} n + O(1).
\end{eqnarray}
Combining (\ref{eq:Pab}) and (\ref{eq:Nstar}) we'd get
\begin{eqnarray*}
  |S| \leq \left( \frac{77}{125} + \frac{583}{1000} \right) n + O(1) = \left( \frac{6}{5} - \frac{1}{1000} \right) n + O(1),
\end{eqnarray*}
a contradiction for $n \gg 0$. It thus suffices to establish (\ref{eq:Pab}). We
divide further into three subcases according to the value of
$x^{*} = \min (P \cap [\alpha n, \, \beta n])$.
\\
\\
\emph{Subcase 1.1:} $x^{*} \leq \frac{n}{5}$. Then $5x^{*}-n \leq 0$ so Lemma
\ref{lem:refred} gives
\begin{eqnarray*}
  |P| \leq n - \frac{5x^{*}}{2} + \frac{1}{2} \leq n - \frac{5\alpha n}{2} + \frac{1}{2} = \frac{29}{50} n + \frac{1}{2} < \frac{583}{1000} n + O(1).
\end{eqnarray*}
If $x^{*} > \frac{n}{5}$, then we can still use Lemma \ref{lem:refred} to get
\begin{eqnarray}\label{eq:xstar}
  |P| \leq |P \cap \left( \frac{2n}{25}, \, \alpha n \right)| + |P \cap [x^{*}, \, 5x^{*} - n)| + n - \frac{5x^{*}}{2} + \frac{1}{2}.
  \end{eqnarray}
  \emph{Subcase 1.2:} $x^{*} > \frac{n}{4}$. The middle interval in
  (\ref{eq:xstar}) is non-empty and we trivially get the bound
  \begin{eqnarray*}
    |P| \leq \left( \alpha - \frac{2}{25} \right)n + (4x^{*} - n) + n - \frac{5x^{*}}{2} + \frac{1}{2} = \frac{11}{125} n + \frac{3x^{*}}{2} + O(1) \\
    \leq \left( \frac{11}{125} + \frac{3\beta}{2} \right) n + O(1) = \frac{583}{1000} n +O(1).
  \end{eqnarray*}
  \emph{Subcase 1.3:} $\frac{n}{5} < x^{*} \leq \frac{n}{4}$. The middle
  interval in (\ref{eq:xstar}) is now empty. The third interval could also
  overlap with the first and Lemma \ref{lem:refred} will give the alternative
  bound
  \begin{eqnarray*}
    |P| \leq n - \frac{5x^{*}}{2} + \max \left\{ 0, \, \min \{\alpha n, \, 5x^{*} - n\} - \frac{2n}{25} \right\} + O(1).
    \end{eqnarray*}
    Ignoring the $O(1)$-term, the right-hand side is otherwise maximized when
    $5x^{*} - n = \alpha n$, which is easily checked to yield
    $|P| \leq \frac{504}{1000} n + O(1)$, which suffices. This completes
    Case 1.
    \\
    \\
    {\sc Case 2:} Neither $P$ nor $N$ meets the interval
    $[\alpha n, \, \beta n]$.
    \\
    \par Here we need to apply more extensively than before the strategy of
    trading off elements of $\pm T$ with those in the rest of $[-n,\,n]$. We
    prove three lemmas, all of which assume $k=5$. The first describes
    trade-offs with $\pm Q$ and is an
    improvement on Lemma \ref{lem:Qtrade} under the conditions of Case 2.
    The second lemma describes
    trade-offs for $\pm M$ and holds for \emph{any} $5$-sum-free set.
    The third lemma describes
    trade-offs for $\pm L$, and makes
    use of the assumption that $|x| \geq \frac{2n}{25}$ for all $x \in S$. 

    To simplify notation, for a $5$-sum-free set $S \subseteq [-n,\,n]$ and $X \in \{P_S, \, N_S\}$, denote
    \begin{eqnarray}\label{eq:gqml}
      g_X = |T \cap X^c|, \;\;\, q_X = |Q \cap X|, \;\;\; m_X = |M \cap X|, \;\;\; l_X = |L \cap X|.
      \end{eqnarray}

    \begin{lemma}\label{lem:Q5trade}
      \begin{eqnarray}\label{eq:q5}
        q_X \leq \frac{2}{5} g_X + 2.
      \end{eqnarray}
    \end{lemma}

    \begin{proof}
      First note that, for $k = 5$, $Q = \left( \frac{4n}{25}, \, \frac{2n}{5} \right]$. Set $Q_{{\hbox{lo}}} = \left( \frac{4n}{25}, \, \alpha n \right)$ and
      $Q_{{\hbox{hi}}} = \left( \beta n, \, \frac{2n}{5} \right]$. So Case 2
      assumes that $X \cap Q = (X \cap Q_{{\hbox{lo}}}) \cup (X \cap Q_{{\hbox{hi}}})$. If $x \in X \cap Q$, let $R_x$ be the interval defined in the proof of
      Lemma \ref{lem:Qtrade}. 

      First suppose $x \in X \cap Q_{{\hbox{lo}}}$. Then $5x < \frac{21n}{25}$
      so, by (\ref{eq:Qspan}),
      \begin{eqnarray}\label{eq:Rlo}
        R_x \subseteq \left( \frac{2n}{5}, \, \frac{11n}{25} \right), \;\;\; {\hbox{if $x \in Q_{{\hbox{lo}}}$}}.
      \end{eqnarray}
      Furthermore, in the notation of Lemma \ref{lem:Qtrade},
      $t \leq n + t_0$, 
      and so
      \begin{eqnarray}\label{eq:Qlo}
        |X \cap Q_{{\hbox{lo}}}| 
        \leq \frac{2g_{X,\,{\hbox{lo}}} + 1}{5} + 1, \;\; {\hbox{where $g_{X,\,{\hbox{lo}}} = \left| X^c \cap \left( \frac{2n}{5}, \, \frac{11n}{25} \right) \right|$}}.
      \end{eqnarray}
      
      Secondly, suppose $x \in X \cap Q_{{\hbox{hi}}}$. Then
      $5x > \frac{33n}{20}$ so, by (\ref{eq:Qspan}),
      \begin{eqnarray}\label{eq:Rhi}
        R_x \subseteq \left( \frac{13n}{20}, \, n \right], \;\;\; {\hbox{if $x \in Q_{{\hbox{hi}}}$}}.
      \end{eqnarray}
      Furthermore, in the notation of Lemma \ref{lem:Qtrade}, $t \geq n + t_0$,
      and so
      \begin{eqnarray}\label{eq:Qhi}
        |X \cap Q_{{\hbox{hi}}}| 
        \leq \frac{2g_{X,\,{\hbox{hi}}} - 1}{5} + 1, \;\; {\hbox{where $g_{X,\,{\hbox{hi}}} = \left| X^c \cap \left( \frac{13n}{20}, \, n \right] \right|$}}.
      \end{eqnarray}
      Since the intervals in (\ref{eq:Rlo}) and (\ref{eq:Rhi}) are disjoint,
      \begin{eqnarray*}
        g_X \geq g_{X, \, {\hbox{lo}}} + g_{X,\,{\hbox{hi}}},
      \end{eqnarray*}
      and then (\ref{eq:q5}) follows from (\ref{eq:Qlo}) and (\ref{eq:Qhi}).
      \end{proof}

      Before going further, we recall that the sumset of two sets $A$ and $B$
      in an abelian group is the set $A+B = \{a+b: a \in A, \, b \in B\}$, and
      that
      \begin{eqnarray}\label{eq:Intsumset}
        |A+B| \geq |A| + |B| - 1, \;\;\; {\hbox{for any finite sets $A, \, B \subseteq \mathbb{Z}$}}.
      \end{eqnarray}
      
      \begin{lemma}\label{lem:M5trade}
        \begin{eqnarray}\label{eq:m5}
          m_X \leq \frac{g_X + 2}{3}.
        \end{eqnarray}
      \end{lemma}

      \begin{proof}
        For $k=5$, $M = \left[ \frac{3n}{25}, \, \frac{4n}{25} \right]$. Let
        $Y = X \cap M$. For any $y_1, \, y_2 \in Y$ one has
        $5y_1 - y_2 \in \left[ \frac{11n}{25}, \, \frac{17n}{25} \right] \subseteq T$. But $5y_1 - y_2 \not\in X$ and so $g_X \geq |5Y - Y|$. To simplify, we
        normalize by setting $Y = y_0 + dY'$, where $y_0$ is the smallest
        element of $Y$ and $d$ their greatest common divisor. So
        $m_X = |Y| = |Y'|$ and $g_X \geq |5Y' - Y'|$.

        For any $0 \leq r \leq 4$, let $Y'_{r} = \{y' \in Y' : y' \equiv r \; ({\hbox{mod $5$}}) \}$. By (\ref{eq:Intsumset}), $5Y' - Y'$ contains at least
        $m_X + |Y'_{r}| - 1$ elements congruent to $-r$ (mod $5$). Let $c$ be
        the number of congruence classes (mod $5$) for which $Y'_{r}$ is
        non-empty. Thus, $|5Y' - Y'| \geq (c+1)m_X - c$. This implies
        (\ref{eq:m5}) if $c \geq 2$. But if $m_X = |Y'| \geq 2$, then
        $c \geq 2$ since gcd$(Y') = 1$. This just leaves the case $m_X = 1$,
        for which (\ref{eq:m5}) just asserts that $g_X \geq 1$. But this is
        also true, since $5y-y = 4y \in T \cap X^c$, for any $y \in Y$.
      \end{proof}

      \begin{lemma}\label{lem:L5trade}
        \begin{eqnarray}\label{eq:l5}
          l_P + l_N \leq \frac{1}{5} ( (g_P + g_N) + 18 ).
        \end{eqnarray}
      \end{lemma}

      \begin{proof}
        For $k=5$ we have $L = \left[1,\,\frac{3n}{25}\right)$. Let
        $L' = \left[1, \, \frac{|T| - 1}{5} \right]$ and $l'_X = |X \cap L'|$.
        It's easy to see that $|L| \leq |L'| + 1$ and so $l_X \leq l'_X + 1$.
        So we need to show that
        \begin{eqnarray}\label{eq:l5prime}
          l'_P + l'_N \leq \frac{g_P + g_N}{5} + \frac{8}{5}.
        \end{eqnarray}

        First consider $P$. If $P \cap L'$ is empty, set $j_P = 0$. Otherwise,
        let $s_P$ be the smallest element of $P \cap L'$
        and $j_P := |T| - 5s_P > 0$. There are $j_P$ numbers in the interval
        $[t_0 + 5s_P, \, n]$. For each such number $\xi$, either
        $\xi \not \in P$ or $\xi - 5s_P \not\in N$, by Lemma \ref{lem:elem}(v).
        Note that each such $\xi$ belongs to the largest $j_P$ elements of $T$,
        whereas $\xi - 5s_P$ belongs to the smallest $j_P$ elements of $T$.
        
        With $s_N$ and $j_N$ defined analogously, there is an interval at the
        top of $T$ consisting of $j_N$ numbers, and for each $\zeta$ in this
        interval, either $\zeta \not\in N$ or $\zeta - 5s_N \not\in P$, where
        moreover $\zeta - 5s_N$ belongs to the smallest $j_N$ elements of $T$.

        For $X \in \{P, \, N\}$ we're assuming $s_X \geq \frac{2n}{25}$, which
        implies that $j_X \leq |T| - \frac{2n}{5}$. But then
        $j_P + j_N \leq 2|T| - \frac{4n}{5} \leq |T|$.
        Because of this, the analysis in the two previous paragraphs implies
        that $g_P + g_N \geq j_P + j_N$. Thus, it suffices to show that
        $l'_X \leq \frac{j_X + 4}{5}$, for
        $X \in \{P,\, N\}$. But the latter follows from the definitions of
        $j_X$. We have $s_X = \frac{|T| - j_X}{5}$ and so
        $X \cap L' \subseteq \left[ \frac{|T| - j_X}{5}, \, \frac{|T| - 1}{5} \right]$. But this interval contains at most
        $\frac{j_X - 1}{5} + 1 = \frac{j_X + 4}{5}$ elements.
\end{proof}

We can now complete the analysis of Case 2. Put $t_X = |X \cap T|$, for $X \in \{P, \, N\}$. Then, using Lemmas \ref{lem:Q5trade}--\ref{lem:L5trade}, we get
\begin{eqnarray*}
  |S| = (t_P + t_N) + (q_P + q_N) + (m_P + m_N) + (l_P + l_N) \\
  \leq \left( 2|T| - (g_P + g_N) \right) + \left( \frac{2(g_P + g_N)}{5} + 4 \right) + \left( \frac{g_P + g_N}{3} + \frac{4}{3} \right) + \left( \frac{g_P + g_N}{5} + \frac{18}{5} \right) \\ \leq |S_{5,\,n}| - \frac{g_P + g_N}{15} + O(1),
\end{eqnarray*}
a contradiction unless $|S \Delta S_{5,\,n}| = O(1)$. This completes the proof
of Theorem \ref{thm:A}.
\qed

    \vspace{3mm}

    \section{{\bf Sums of dilates}}\label{sect:dilates}

    A well-known result of Shakan \cite{S} on sums of dilates states that, for any
    vector $\overline{\lambda} = (\lambda_1,\,\dots,\,\lambda_m) \in \mathbb{Z}^m$ such that GCD$(\lambda_1,\,\dots,\,\lambda_m) = 1$, for all finite sets $A$
    of integers one has
    \begin{eqnarray}\label{eq:bukh}
      |\lambda_1 A + \dots + \lambda_m A| \geq \left( \sum_{i=1}^{m} |\lambda_i| \right) \cdot |A|  - O_{\overline{\lambda}}(1).
    \end{eqnarray}
    As a special case, $|4A - A| \geq 5|A| - O(1)$. Note that this is a much
    stronger inequality than that obtained from the basic sumset bound
    (\ref{eq:Intsumset}) since, for a set $A$ of integers, $|4A| = |A|$ and so
    (\ref{eq:Intsumset}) would just give $|4A - A| \geq 2|A| - 1$.

    For the proof of Theorem \ref{thm:B} we need a ``one-sided'' version of
    Shakan's inequality, which gives a non-trivial lower bound on the number of
    non-negative elements in $4A - A$, when $A$ is also a set of non-negative
    integers.
    \\
    \par
    Define the function $\psi: \mathbb{R}_{+} \rightarrow \mathbb{R}_{+}$ by
    \begin{eqnarray}\label{eq:g}
      \psi(s) = 60 (s+1)^{4/5} - 56.
    \end{eqnarray}
    This function has the following properties:
    \par (i) $\psi(0) = 4$, $\psi(s)$ is an increasing function of $s$ and $\psi(s) = O(s^{4/5}) = o(s)$.
    \par (ii) For all $s \geq 1$, $\psi(s) > 3 \cdot \psi(s/4)$.
    \par (iii) For all $s \geq 3$, $\psi(s) - \psi(3s/4) > 3$.
    \\
    \\
    All three properties are easily verified. Property (i) is trivial. For
    (ii), note that it is equivalent to having $h(s) < \frac{28}{15}$ for all
    $s \geq 1$, where $h(s) = 3 \left( \frac{s}{4} + 1 \right)^{4/5} - (s+1)^{4/5}$. One can check directly that $h(1) < \frac{28}{15}$ and that
    $h^{\prime}(s) < 0$ for $s > 0$. Similarly, for (iii) we set
    $k(s) = \psi(s) - \psi(3s/4)$, verify directly that $k(3) > 3$ and then
    check that $k^{\prime}(s) > 0$ for all $s > 0$.

    For a finite set $A$ of integers and a real number $r$, set
    \begin{eqnarray*}
      P(A, \, r) = |(4A - A) \cap (r, \, \infty)| = \#\{w \in 4A - A: w > r\}.
    \end{eqnarray*}
    Put $P(A) = P(A,\,0)$ and $d(A) = 4|A| - P(A)$. So $P(A)$ is the number of
    strictly positive elements in $4A - A$. 
    
    \begin{theorem}\label{thm:dilates}
      Let $A$ be a finite set of non-negative integers, with $0 \in A$. If $s$
      is the largest element of $A$, then $P(A) \geq 4|A| - \psi(s)$,
      equivalently $d(A) \leq \psi(s)$.
    \end{theorem}

    \begin{remark}\label{rem:dilates}
      The theorem says that, for any finite set $A$ of non-negative integers
      containing zero, the set $4A - A$ contains at least $4|A| - O(s^{4/5})$
      positive elements. This set also contains at least $|A|$ non-positive
      elements, namely $0-a$ for all $a \in A$. Hence
      $|4A - A| \geq 5|A| - O(s^{4/5})$. This is weaker than in Shakan's result.
      We would recover the $\overline{\lambda}=(4,\,-1)$ case of the latter if
      Theorem \ref{thm:dilates} held with
      $d(A) = O(1)$. This remains open. Note that, if we
      could take $d(A) = O(1)$, the analysis
      from this point onwards would yield Conjecture \ref{conj:C}.
    \end{remark}

    \emph{Proof of Theorem \ref{thm:dilates}.} The proof is by strong
    induction on
    $s$. If $s=0$, then $A = \{0\}$, $|A| = 1$, $P(A) = 0$ and
    $d(A) = 4 = \psi(0)$. Now suppose $s \geq 1$.
    \\
    \par For each $r \in \{0,\,1,\,2,\,3\}$, let $A_r = A \cap (r + 4\mathbb{Z})$. The sets $4A_r - A_r$ are pairwise disjoint subsets of $4A - A$, since they
    lie in distinct residue classes mod $4$. Thus, $P(A) = \sum_{r} P(A_r)$. We
    normalize by writing $A_r = 4B_r + r$, $\mu_r = \min B_r$ and
    $\hat{B}_r = B_r - \mu_r$. It is easily checked that
    $w \mapsto 4w + 12\mu_r + 3r$ maps the positive part of
    $4\hat{B}_r - \hat{B}_r$ injectively into that of $4A_r - A_r$. Hence,
    $P(A_r) \geq P(\hat{B}_r)$ and so $P(A) \geq \sum_{r} P(\hat{B}_r)$. Since
    $|\hat{B}_r| = |A_r|$ and $\sum_r |A_r| = |A|$, this can be rewritten as
    \begin{eqnarray}\label{eq:dhatr}
      d(A) \leq \sum_{r} d(\hat{B}_r).
    \end{eqnarray}
    Note that the largest element of any $\hat{B}_r$ is at most $s/4$, which is
    strictly less than $s$ whenever $s \geq 1$. If some $A_r$ is empty, then
    Property (ii) of the function $\psi(s)$ suffices to complete the induction
    step.

    So we may assume $A$ meets every residue class mod $4$. In particular, this
    means $s \geq 3$. In this case, for any real $\theta \in \left[0, \, \frac{3s}{4} \right)$, we claim that
    \begin{eqnarray}\label{eq:windows}
      d(A) \leq d(A_{\leq \theta}) + 4\left| A \cap \left( \theta, \, \theta + \frac{s}{4} \right] \right| - |A| + 4, \;\;\; {\hbox{where $A_{\leq \theta} = A \cap [0,\,\theta]$}}.
    \end{eqnarray}
    We show this by considering different, almost non-overlapping parts of
    $4A - A$. Firstly, the positive part of $4A_{\leq \theta} - A_{\leq \theta}$
    lies in $(0,\,4\theta]$. Secondly, for each residue class $r$, let
    $s_r = \max A_r$ and $E_r = E_{r,\,1} \cup E_{r,\,2}$, where
    \begin{eqnarray*}
      E_{r,\,1} = \left\{ 4a - s_r : a \in A, \, a > \theta + \frac{s_r}{4} \right\}, \;\;\; E_{r,\,2} = \{ 4s - a: a \in A_r \}.
    \end{eqnarray*}
    It is easily checked that
    \begin{itemize}
    \item
      the two parts of $E_r$ share only the single element $4s - s_r$,
    \item
      every element of $E_r$ is congruent to $-r \; ({\hbox{mod $4$}})$,
    \item
      every element of $E_r$ is greater than $4\theta$.
    \end{itemize}
    Thus,
    \begin{eqnarray}\label{eq:Pwindows}
      P(A) \geq P(A_{\leq \theta}) + \sum_{r=0}^{3} (|E_{r,\,1}| + |E_{r,\,2}| - 1) \nonumber \\ = P(A_{\leq \theta}) + \sum_{r=0}^{3} \left| A \cap \left( \theta + \frac{s_r}{4}, \, s \right] \right| + |A_r| - 1 \nonumber \\ \Rightarrow \; P(A) \geq P(A_{\leq \theta}) + 4 \left| A \cap \left( \theta + \frac{s}{4}, \, s \right] \right| + |A| - 4. 
    \end{eqnarray}
    Upon noting that 
    \begin{eqnarray*}
      |A| = |A_{\leq \theta}| + \left| A \cap \left( \theta, \, \theta + \frac{s}{4} \right] \right| + \left| A \cap \left( \theta + \frac{s}{4}, \, s \right] \right|,
    \end{eqnarray*}
    it is easily verified that (\ref{eq:Pwindows}) can be rewritten as
    (\ref{eq:windows}).

    Now the four windows $\left( \theta, \, \theta + \frac{s}{4} \right]$, for
    $\theta \in \left\{0, \, \frac{s}{4}, \, \frac{s}{2}, \, \frac{3s}{4} - \frac{1}{2} \right\}$, together cover all of $A$ except for the elements $0$ and
    $s$. They are almost pairwise disjoint, except that at most one integer may
    appear in the last two windows. Hence, some window must contain at most
    $\frac{|A| - 1}{4}$ elements. For this choice of $\theta$,
    (\ref{eq:windows}) gives
    \begin{eqnarray*}
      d(A) \leq d(A_{\leq \theta}) + 3, \;\;\; {\hbox{for some $\theta < \frac{3s}{4}$}}.
    \end{eqnarray*}
    Now Property (iii) of the function $\psi(s)$ suffices to complete the
    induction step, in the case where $A$ meets every residue class mod $4$.
    \qed

    \begin{corollary}\label{cor:dilates}
      Let $A$ be a finite set of non-negative integers, with $0 \in A$ and $s = \max A$. Then, for any integer $K$, 
      \par (i) $P(A,\,K) \geq 4|A| - \psi(s) - \max \{K,\,0\}.$
      \par (ii) If $A$ meets exactly $c$ residue classes mod $4$, then
      \begin{eqnarray}\label{eq:dilatesK}
        P(A,\,K) \geq 4|A| - \frac{c}{4} \cdot \max\{K,\,0\} - c \cdot \psi \left( \frac{s}{4} \right).
      \end{eqnarray}
    \end{corollary}

    \begin{proof}
      Part (i) is a trivial consequence of Theorem \ref{thm:dilates}, since
      increasing
      the threshold by one can remove at most one element of $4A - A$. For part
      (ii), consider the map $w \mapsto 4w + 12\mu_r + 3r$ in the first part of
      the proof of the theorem. The image exceeds $K$ as soon as
      $w > \frac{K}{4}$, so
      \begin{eqnarray*}
        P(A_r,\,K) \geq P\left( \hat{B}_r,\, \lfloor \frac{K}{4} \rfloor \right) \geq 4|A_r| - \psi \left( \frac{s}{4} \right) - \frac{1}{4} \cdot \max \{K,\,0\},
      \end{eqnarray*}
      where the second inequality uses part (i), the fact that $\psi (\cdot)$
      is an increasing function and that $|A_r| = |\hat{B}_r|$. Summing over
      the $c$ classes met by $A$ yields (\ref{eq:dilatesK}), since the
      contributions to $P(A,\,K)$ come from distinct classes mod $4$.
    \end{proof}
    
    \vspace{3mm}

    \section{{\bf Case of $k=4$: Proof that $S_4$ achieves maximum density}}\label{sect:kfour}

In this section we will prove that, if $S$ is a $4$-sum-free subset of
$[-n,\,n]$, then
\begin{eqnarray}\label{eq:Off}
  |S| \leq \frac{33}{32} n + o(n).
\end{eqnarray}
For $k=4$, we have
\begin{eqnarray*}
  L = \left[1,\, \frac{n}{8} \right), \;\;\; M = \left[ \frac{n}{8}, \, \frac{n}{4} \right], \;\;\; Q = \left( \frac{n}{4}, \, \frac{n}{2} \right], \;\;\; T = \left( \frac{n}{2}, \, n \right].
\end{eqnarray*}
The basic strategy is still to show that trading elements of $T$ for those of
$L, \, M$ or $Q$ always goes badly in a $4$-sum-free set. Lemmas
\ref{lem:L4trade}, \ref{lem:M4trade} and \ref{lem:Q4trade} make this strategy
precise, and are
analogues of Lemmas \ref{lem:Q5trade}--\ref{lem:L5trade}, though in a different
order. The argument for $Q$ is much more intricate now and that for $M$ relies
on Theorem \ref{thm:dilates}. The argument for $L$ is the simplest of the three
in the current setting. 

Having established our three trade-off lemmas, plus Lemma \ref{lem:4refred}
which is a simple modification of Lemma \ref{lem:refred}, we divide the proof
of (\ref{eq:Off}) into four cases, depending on the smallest elements of
$P \cap M$ and $N \cap M$. The aforementioned lemmas allow us to deal quickly
with Cases 1--3. Case 4, however, requires us to dissect the intervals $M$ and
$Q$ further, especially $M$, and analyze trade-offs for the different
parts in much greater
detail. The analysis yields seven further inequalities, Lemmas
\ref{lem:R}--\ref{lem:U1}, the proofs of several of which rely on Theorem \ref{thm:dilates}. The proof of (\ref{eq:Off}) is then completed by
combining all our inequalities appropriately: this last step is conceptually
simple, but computationally intimidating. Note that while Case 4 requires by
far the most intricate analysis, extremal $4$-sum-free sets never arise there:
see (\ref{eq:C4bound}) below. Nor can they arise in Case 3. Thus, when we come
to complete the proof of
Theorem \ref{thm:B} in Section \ref{sect:stability}, only Cases 1 and 2 will be
in focus. 

\subsection*{Four preparatory lemmas} As before, $P$ and $N$ denote the positive and negative parts of a
$4$-sum-free set $S$, i.e.: $S = P \cup (-N)$. For $X \in \{P, \, N\}$ we adopt
the same notation as in (\ref{eq:gqml}) and also $\tau_X := |X \cap T|$, so that
$g_X = |T| - \tau_X$. For any set of integers $Y$, we denote
\begin{itemize}
\item
  $Y_r := Y \cap (r + 4\mathbb{Z})$, for $0 \leq r \leq 3$,
\item
  $Y^{\prime}_{r} := Y \cap (r + 8\mathbb{Z})$, for $0 \leq r \leq 7$.
\end{itemize}

Set $\tau_r = \tau_{r,\,X} := |(X \cap T)_r|$ and $\tau_e := \tau_0 + \tau_2$. Finally, set $\tau^{\prime}_{r} := |(X \cap T)^{\prime}_{r}|$.
\\
\par We begin with the trade-off lemma for $L$, which is the simplest of the
three. 

\begin{lemma}\label{lem:L4trade}
At least one of the following must hold:
  \\ (a) \begin{eqnarray}\label{eq:l4a}
             l_P + l_N \leq \frac{1}{4} (g_P + g_N) + O(1).
           \end{eqnarray}
           (b) \begin{eqnarray}\label{eq:l4b}
                 \tau_P + \tau_N \leq \frac{n}{2} + O(1).
               \end{eqnarray}
\end{lemma}

\begin{proof}
  By definition,
               \begin{eqnarray}\label{eq:lhelp}
                 \tau_P + \tau_N = \sum_{j=0}^{3} (\tau_{j,\,P} + \tau_{j,\,N}).
               \end{eqnarray}
               Pick the congruence class $j \; ({\hbox{mod $4$}})$ for which
               the $j$:th right-hand term is largest, and consider two cases.
               \par {\sc Case 1:} Either $(P \cap T)_j$ or $(N \cap T)_j$ is
               empty. Now, $|T| \leq \frac{n}{2} + O(1)$, so
               $\max \{\tau_{j,\,P}, \, \tau_{j,\,N} \} \leq \frac{n}{8} + O(1)$.
               So in Case 1, each term of the sum is at most
               $\frac{n}{8} + O(1)$, and thus (b) holds.
               \par {\sc Case 2:} Both $(P \cap T)_j$ and $(N \cap T)_j$ are
               non-empty sets. Let $D$ denote the difference set
               $(P \cap T)_j - (N \cap T)_j$. By (\ref{eq:Intsumset}) it
               contains at least $\tau_{j,\,P} + \tau_{j,\,N} - 1 \geq \frac{\tau_P + \tau_N}{4} - 1$ elements. By Lemma \ref{lem:elem}(v), $D$ is disjoint from
               both $4(P \cap L)$ and $-4(N \cap L)$. Hence, the three sets
               $D$, $4(P \cap L)$ and $-4(N \cap L)$ are pairwise disjoint, and
               all consist of multiples of $4$ in the interval
               $\left( -\frac{n}{2}, \, \frac{n}{2} \right)$. This interval
               contains $\frac{n}{4} + O(1)$ multiples of $4$ and thus,
               \begin{eqnarray*}
                 l_P + l_N + \frac{\tau_P + \tau_N}{4} - 1 \leq \frac{n}{4} + O(1).
               \end{eqnarray*}
               Since $\tau_X = \frac{n}{2} - g_X + O(1)$, for $X \in \{P,\,N\}$,
               this implies (a).
\end{proof}

Next we turn to trade-offs for $M$, where our basic result will follow quickly
from Theorem \ref{thm:dilates}.

\begin{lemma}\label{lem:M4trade}
  Let $X \in \{P,\,N\}$ and suppose $X \cap M \neq \phi$. Let $u$ be the
  smallest element of $X \cap M$. Then
  \begin{eqnarray}\label{eq:M4trade}
    m_X \leq \frac{g_X}{4} + \kappa(u) + O(n^{4/5}), \;\;\; {\hbox{where $\kappa(u) = \max \left\{ 0, \, \frac{n/2 - 3u}{4} \right\}$}}.
  \end{eqnarray}
\end{lemma}

\begin{proof}
  Let $A = (X \cap M)$ and $A^{*} = A - u$, so $A^{*} \subseteq \mathbb{Z}_{\geq 0}$, $0 \in A^{*}$, $\max A^{*} \leq n/8$ and $|A^{*}| = m_X$. The set $4A - A$ is disjoint from
  $X$, by Lemma \ref{lem:elem}(ii). Every number in this set is less than $n$
  and,
  since $4y - z = 4(y-u) - (z-u) + 3u$, the number of its elements exceeding
  $n/2$ is $P(A^{*},\,K)$, where $K = \lfloor \frac{n}{2} \rfloor - 3u$. Hence, by
  Corollary \ref{cor:dilates}(i),
  \begin{eqnarray*}
    g_X \geq P(A^{*},\,K) \geq 4m_X - \max \{K,\,0\} - \psi \left( \frac{n}{8} \right).
  \end{eqnarray*}
  Since $\frac{1}{4} \cdot \max \{K,\,0\} \leq \kappa(u)$ and $\psi \left(\frac{n}{8}\right) = O(n^{4/5})$, this gives (\ref{eq:M4trade}). 
\end{proof}

\begin{remark}\label{rem:M4trade}
  Since $u \geq n/8$, $\kappa(u) \leq n/32$, which already hints at why the set
  $S_4$ achieves maximum density. 
\end{remark}

We turn to trade-offs between $T$ and $Q$, where the argument is now quite
intricate. We begin with a definition:

\begin{definition}\label{defi:exceptional}
  Let $X \in \{P, \, N\}$, $A = X \cap T$ and $B = X \cap Q$. The set $X$ is
  said to be \emph{exceptional} if exactly one of $A_1$ and $A_3$ is an empty set and
  \begin{itemize}
  \item
    If $A_3 = \phi$, then also $B_3 = \phi$ and $\tau^{\prime}_{6} \geq \tau_0 + \tau_{2}^{\prime} + 1$.
  \item
    If $A_1 = \phi$, then also $B_1 = \phi$ and $\tau_{2}^{\prime} \geq \tau_0 + \tau^{\prime}_{6} + 1$.
  \end{itemize}
\end{definition}

\begin{lemma}\label{lem:Q4trade}
  For $X \in \{P, \, N\}$,
  \par (a) If $X$ is not exceptional, then
  \begin{eqnarray}\label{eq:q4}
    q_X \leq \frac{g_X}{2} + O(1).
  \end{eqnarray}
  \par (b) If $X$ is exceptional, then
  \begin{eqnarray}\label{eq:exceptineq}
    |X| \leq \frac{393}{880} n + o(n).
  \end{eqnarray}
\end{lemma}

\begin{proof}
  Let $C = Q \backslash B = Q \cap X^c$. Since $|Q| < \frac{|T|}{2} + 1$, any
  bound of the form
  \begin{eqnarray}\label{eq:q4c}
    |C| \geq \frac{\tau_X}{2} - O(1)
  \end{eqnarray}
  will imply (\ref{eq:q4}). We will make repeated use of the following two
  facts:
  \\
  \\
  \emph{(F1)}: Let $W = \{t/2 : t \in A_0 \cup A_2\}$. Then $W \subseteq C$, since $4\left( \frac{t}{2} \right) = t +t$. Here $|W| = \tau_e$ and $W = W_{e} \cup W_{o}$,
  where $W_e = A_0/2$ and $W_o = A_2/2 = A^{\prime}_{2}/2 \cup A^{\prime}_{6}/2$.
  \\
  \emph{(F2)}: If $t, \, t' \in A$ and $4 \, | \, (t+t')$, then $\frac{t+t'}{4} \in C$.
  \\
  \par First suppose $X$ is not exceptional. By (F2), each of $A_0 + A_0$,
  $A_2 + A_2$ and $A_1 + A_3$ is contained in the dilate $4C$. These sets
  need not be disjoint, but using
  (\ref{eq:Intsumset}) and taking a maximum, it easily follows that
  (\ref{eq:q4c}) holds if any of
  the following apply:
  \begin{itemize}
  \item
    $A_1$ and $A_3$ are both non-empty,
  \item
    $A_1$ and $A_3$ are both empty (in which case $\tau_X = \tau_0 + \tau_2$),
  \item
    $\max \{ \tau_0, \, \tau_2 \} \geq \frac{\tau_X}{4}$.
  \end{itemize}

  We can thus assume $\tau_0$ and $\tau_2$ are both less than
  $\frac{\tau_X}{4}$, and that exactly one of $A_1$ and $A_3$ is non-empty. We
  assume $A_3$ is empty and observe that the rest of the argument can be
  immediately reformulated if instead $A_1$ is empty. We divide the
  analysis into two cases, depending on whether $B_3$ is empty or not. 
  \\
  \\
  {\sc Case 1:} $B_3$ is empty. Since $X$ is not exceptional, we must have
  $\tau_0 + \tau_{2}^{\prime} \geq \tau^{\prime}_{6}$, which is the same as saying
  that $\tau_0 + \tau^{\prime}_{2} \geq \frac{\tau_e}{2}$. Since $B_3$ is empty,
  we have $Q_3 \subseteq C$. As noted in (F1), $W_e$ and $A^{\prime}_{2}/2$ are
  also in $C$, and both are disjoint from $Q_3$ as well as from one another.
  Thus,
  \begin{eqnarray*}
    |C| \geq |Q_3| + |W_e| + |A^{\prime}_{2}| \geq \left( \frac{|T|}{8} - O(1) \right)  + \tau_0 + \tau^{\prime}_{2} \geq \frac{\tau_1 + \tau_e}{2} - O(1) = \frac{\tau_X}{2} - O(1),
  \end{eqnarray*}
  where the last equality is due to $A_3$ being empty. This establishes
  (\ref{eq:q4c}).
  \\
  \\
  {\sc Case 2:} $B_3$ is non-empty. Let $b = \max B_3$ and set
  \begin{eqnarray*}
    Y = Q_3 \cap \left(a, \, \frac{n}{2} \right], \;\;\; Z = \{ \frac{b+a}{4} : a \in A_1, \, a > n-b\}.
  \end{eqnarray*}
  $Y$ is a subset of $C$ by definition of $b$, and so is $Z$ since
  $Z \subseteq \left( \frac{n}{4}, \, \frac{n+b}{4} \right]$. The elements of
  $A_1$ below $n-b$ lie in a single class mod 4 of the interval
  $\left( \frac{n}{2}, \, n-b\right]$, which has the same length as
  $\left(b, \, \frac{n}{2} \right]$. Thus, $|Y| + |Z| \geq \tau_1 - 1$. Let
  $i = |Y \cap Z|$ and $\omega = |W \cap (Y \cup Z)|$. Note that
  $i \leq i^{*} := \left| Q_3 \cap \left(b, \, \frac{n+b}{4} \right] \right|$.
  Thus, 
  \begin{eqnarray*}
    |C| \geq |W \cup Y \cup Z| \geq |W| + (|Y| + |Z|) - i - \omega \geq \tau_e + (\tau_1 - 1) - i^{*} - \omega = \tau_X - 1 - i^{*} - \omega,
  \end{eqnarray*}
  where the last equality is due to $A_3$ being empty. Thus (\ref{eq:q4c})
  holds if $\omega \leq \frac{\tau_X}{2} - i^{*} + O(1)$. So we may assume
  otherwise. 

  Next, $Y$ consists of odd numbers, so
  \begin{eqnarray*}
    \omega \leq |W_e \cap Z| + |W_o \cap (Y \cup Z)| \leq |W_e \cap Z| +|W_o| = |W_e \cap Z| + \tau_2 < |W_e \cap Z| + \frac{\tau_X}{4}.
  \end{eqnarray*}
  So we are assuming that $|W_e \cap Z| > \frac{\tau_X}{4} - i^{*} - O(1)$. Let
  \begin{eqnarray*}
    U = \left\{ \frac{a}{4} : a \in A_0, \, \frac{a}{2} \in W_e \cap \mathbb{Z} \right\}.
  \end{eqnarray*}
  Then $|U| = |W_e \cap Z|$ and $U \subseteq \left( \frac{n}{8}, \, \frac{n+b}{8} \right]$. By (F2), $U+U \subseteq C$. But $\max (U+U) \leq \frac{n+b}{4}$, so $|(U+U) \cap Y| \leq i^{*}$. Thus, and using also (\ref{eq:Intsumset}),
  \begin{eqnarray*}
    |C| \geq |(U+U) \cup Y| \geq 2|U| - 1 + |Y| - i^{*} = 2|W_e \cap Z| - 1 + |Y| - i^{*} \\ > 2\left( \frac{\tau_X}{4} - i^{*} - O(1) \right) - 1 + |Y| - i^{*} = \frac{\tau_X}{2} - O(1) - (3i^{*} - |Y|).
  \end{eqnarray*}
  So it suffices to show an $O(1)$ upper bound on $3i^{*} - |Y|$. But recall
  that $i^{*} = \left| Q_3 \cap \left( b, \, \frac{n+b}{4} \right] \right|$. If
  $b \geq \frac{n}{3}$, then $i^{*} = 0$ and $3i^{*} - |Y| \leq 0$. If
  $b < \frac{n}{3}$, then $i^{*} \leq 1 + \frac{1}{4} \left( \frac{n+b}{4} - b \right) = 1 + \frac{n-3b}{16}$. Finally, recall that
  $Y = Q_3 \cap \left( b, \, \frac{n}{2} \right]$, so
  $|Y| \geq \frac{1}{4} \left( \frac{n}{2} - b \right) - 1 = \frac{n-2b}{8} - 1$. Thus, assuming $b < \frac{n}{3}$,
  \begin{eqnarray*}
    3i^{*} - |Y| \leq 3\left( \frac{n-3b}{16} + 1\right) - \left( \frac{n-2b}{8} - 1 \right) = 4 + \frac{n-5b}{16} < 4,
  \end{eqnarray*}
  since $b > \frac{n}{4} > \frac{n}{5}$. This completes the proof of part (a)
  of Lemma \ref{lem:Q4trade}. 
  \\
  \par Now suppose $X$ is exceptional. We present the analysis in the case when
  $A_3 = B_3 = \phi$ and $\tau^{\prime}_{6} \geq \tau_0 + \tau^{\prime}_{2} + 1$,
  leaving the reader to check that a directly analogous argument applies in the
  other exceptional case. We will show that, under these conditions,
  \begin{eqnarray}\label{eq:over8th}
    \left| X \cap \left( \frac{n}{8}, \, n \right] \right| \leq \frac{3n}{8} + O(1).
  \end{eqnarray}
  In \cite{H1} it was proven that a $4$-sum-free subset of $[1,\,n]$ has size
  at most $\frac{441}{770} n + o(n)$. Thus (\ref{eq:over8th}) would imply that
  \begin{eqnarray*}
    |X| \leq \frac{1}{8} \left( \frac{441}{770} n + o(n) \right) + \frac{3}{8} n + O(1) = \frac{393}{880} n + o(n),
  \end{eqnarray*}
  as desired. So it just remains to establish (\ref{eq:over8th}).

  For each $y \in Q$, define the \emph{chain} $C(y)$ as
  \begin{eqnarray*}
    C(y) = \left\{ \begin{array}{lr} \{y,\,2y\}, & {\hbox{if $y$ is odd}}, \\ \left\{ \frac{y}{2}, \, y, \, 2y \right\}, & {\hbox{if $y$ is even}}. \end{array} \right.
  \end{eqnarray*}
  The chains partition $\left( \frac{n}{8}, \, n \right] \backslash (T_1 \cup T_3)$. Since $X$ never contains both $x$ and $2x$, $|X \cap C(y)| \leq 1$, except
  when $y$ is even and both $y/2$ and $2y$ lie in $X$. Let
  $K_0 = \{x \in \left( \frac{n}{8}, \, \frac{n}{4} \right]: 4x \in X\}$ and
  note that $|K_0| \leq \tau_0$.

  If $y \in Q_3$, then $y$ is odd and $|X \cap C(y)| = 1$ exactly when
  $2y \in A^{\prime}_6$, so the number of such chains is exactly
  $\tau^{\prime}_6$. Let $\varepsilon$ be the number of empty chains $C(y)$, for
  $y \in Q_0 \cup Q_1 \cup Q_2$. Since $A_3$ is empty, the above analysis shows
  that
  \begin{eqnarray*}
    \left| X \cap \left( \frac{n}{8}, \, n \right] \right| \leq \tau_1 + \tau^{\prime}_{6} + |Q_0 \cup Q_1 \cup Q_2| - \varepsilon + |K_0| \leq \left( \frac{n}{8} + O(1) \right) + \tau^{\prime}_{6} + \left( \frac{3n}{16} + O(1) \right) - \varepsilon + \tau_0,
  \end{eqnarray*}
  so it suffices now to show that
  \begin{eqnarray}\label{eq:taucount}
    \tau^{\prime}_{6} + \tau_0 - \varepsilon \leq \frac{n}{16} + O(1).
  \end{eqnarray}
  Let $V$ be the set such that $2(4V + 3) = A^{\prime}_{6}$. In other words,
  $a \mapsto \frac{1}{4} \left( \frac{a}{2} - 3 \right)$ is a bijection from
  $A^{\prime}_{6}$ to $V$, so $|V| = |A^{\prime}_{6}| = \tau^{\prime}_{6}$.

  If $V$ misses a parity class (i.e.: a residue class mod $2$), then
  $A^{\prime}_{6}$ will lie in a single residue class mod $16$, so
  $\tau^{\prime}_{6} \leq \frac{n}{32} + O(1)$. But the assumptions of
  exceptionality also imply that $\tau_0 < \tau^{\prime}_{6}$, so
  (\ref{eq:taucount}) follows in this case.

  Finally, suppose $V$ contains both even and odd numbers, and let $V_e$ and
  $V_o$ denote the respective non-empty sets. For $v \in V_e$ and
  $v' \in V_o$, let $y = y(v,\,v') = \frac{1}{2} ((3+4v) + (3+4v'))$ and note
  that
  \begin{itemize}
  \item
    $y \equiv 1 \; ({\hbox{mod $4$}})$,
  \item
    $y \in Q$,
  \item
    $y \not\in X$, by (F2).
  \end{itemize}
  Thus, $C(y) = \{y,\,2y\}$ and $|X \cap C(y)| = 1$ if and only if
  $2y \in A^{\prime}_{2}$. What this shows is that
  $\varepsilon \geq |V_e + V_o| - \tau^{\prime}_{2}$. But $|V_e + V_o| \geq |V_e| + |V_o| - 1 = |V| - 1 = \tau^{\prime}_{6} - 1$, by (\ref{eq:Intsumset}), and so
  $\varepsilon \geq \tau^{\prime}_{6} - 1 - \tau^{\prime}_{2}$. Thus,
  $\tau^{\prime}_{6} + \tau_0 - \varepsilon \leq \tau_0 + \tau^{\prime}_{2} + 1$,
  which in turn is at most $\tau^{\prime}_{6}$, by the assumptions of
  exceptionality. But $\tau^{\prime}_{6} \leq \frac{n}{16} + O(1)$ trivially,
  and this establishes (\ref{eq:taucount}).

  The proof of Lemma \ref{lem:Q4trade} is now complete.
\end{proof}

The final preparatory lemma is in the same spirit as Lemma \ref{lem:refred}.

\begin{lemma}\label{lem:4refred}
  If $x \in P \cap M$, then $\tau_P + |N \cap [1,\,n-4x]| \leq |T|$. The same
  holds with $P$ and $N$ interchanged.
\end{lemma}

\begin{proof}
  Note that $\frac{n}{2} \leq 4x \leq n$.
  If $y \in P \cap T$ then $4x-y \not\in S$, by Lemma
  \ref{lem:elem}(ii). If $y > 4x$, this means $y-4x \not\in N$. So
  $y \mapsto y-4x$ maps $P \cap (4x,\,n]$ injectively into
  $[1,\,n-4x] \backslash N$. Thus,
  \begin{eqnarray*}
    |P \cap (4x,\,n]| + |N \cap [1,\,n-4x]| \leq n-4x.
  \end{eqnarray*}
  Adding $\left| P \cap \left( \frac{n}{2}, \, 4x \right] \right|$ to both
  sides yields the result.
  \end{proof}

  \subsection*{Proof of (\ref{eq:Off}), reduction to four cases} Let $F(n)$ be
  the maximum size of a $4$-sum-free subset of $[-n,\,n]$. Suppose we can show
  that there is a non-decreasing function $\varepsilon(n)$ such that
  $\varepsilon(n) = o(n)$ and, for any $4$-sum-free subset
  $S \subseteq [-n,\,n]$, either
  \par (E1): $|S| \leq \frac{33}{32} n + \varepsilon(n)$, or
  \par (E2): $|S| \leq \frac{n}{2} + C + F\left( \lfloor \frac{n}{2} \rfloor \right)$, for some constant $C$.
  \\
  \\
  We claim this will suffice. Put $E(n) := F(n) - \frac{33}{32} n$. Then we
  would have
  \begin{eqnarray*}
    E(n) \leq \max \left\{ \varepsilon(n), \, E\left( \lfloor \frac{n}{2} \rfloor \right) + C - \frac{n}{64} \right\}.
  \end{eqnarray*}
  A simple induction will then yield $E(n) \leq \max \{\varepsilon(n), \, C'\}$,
  where $C' = \max_{n < 64C} E(n)$.
  \\
  \par Let $S$ be a $4$-sum-free subset of $[-n,\,n]$ and first suppose that
  either $P$ or $N$ is exceptional. Applying Lemma \ref{lem:Q4trade}(b) to the
  appropriate side, and the bound of \cite{H1} to the other side, we get
  \begin{eqnarray}\label{eq:excbd}
    |S| \leq \left( \frac{393}{880} + \frac{441}{770} \right) n + o(n) = \left( \frac{33}{32} - \frac{147}{12320} \right) n + o(n),
  \end{eqnarray}
  which yields (E1). Secondly, suppose Lemma \ref{lem:L4trade}(b) holds, with a
  constant $C$ in the $O(1)$-term. Then,
  \begin{eqnarray*}
    |S| = \tau_P + \tau_N + \left| S \cap \left[ -\frac{n}{2}, \, \frac{n}{2} \right] \right| \leq \frac{n}{2} + C + F\left( \lfloor \frac{n}{2} \rfloor \right),
  \end{eqnarray*}
  which is (E2). So we can henceforth assume Lemma \ref{lem:L4trade}(a) holds
  and that Lemma \ref{lem:Q4trade}(a) holds for both $P$ and $N$. Thus,
  \begin{eqnarray*}
    |S| = (\tau_P + \tau_N) + (q_P + q_N) + (m_P + m_N) + (l_P + l_N) \\
    \leq (2|T| - (g_P + g_N)) + \frac{1}{2} (g_P + g_N) + (m_P + m_N) + \frac{1}{4} (g_P + g_N) + O(1) \\
    \leq n + (m_P + m_N) - \frac{1}{4} (g_P + g_N) + O(1),
  \end{eqnarray*}
  so for (E1) to hold it suffices to have
  \begin{eqnarray}\label{eq:E1bound}
    m_P + m_N \leq \frac{1}{4} (g_P + g_N) + \frac{n}{32} + o(n).
  \end{eqnarray}
  Now we can distinguish our four cases.
  \\
  \\
  {\sc Case 1:} Either $P \cap M$ or $N \cap M$ is empty. Then (\ref{eq:E1bound}) follows immediately from Lemma \ref{lem:M4trade} and Remark \ref{rem:M4trade}.
  \\
  \par In the remaining cases, we may assume both $P \cap M$ and $N \cap M$ are
  non-empty. Let $a$ and $b$ be the smallest elements of $P \cap M$ and
  $N \cap M$, respectively. By definition of $M$,
  $\frac{n}{8} \leq a, \, b \leq \frac{n}{4}$. Lemma \ref{lem:M4trade} tells us
  that
  \begin{eqnarray}\label{eq:E1M}
    m_P + m_N \leq \frac{1}{4} (g_P + g_N) + \kappa(a) + \kappa(b) + o(n).
    \end{eqnarray}
    {\sc Case 2:} $\max \{a,\,b\} > \frac{n}{6}$. WLOG, say $a > \frac{n}{6}$.
    Then $\kappa(a) = 0$ and, since $\kappa(b) \leq \frac{n}{32}$, we're done.
    \\
    \\
    {\sc Case 3:} $a+b \leq \left( (\frac{33}{128} - c \right) n$, for any fixed $c \geq 0$. In
    this case, we can bypass (\ref{eq:E1M}) and just apply Lemma
    \ref{lem:4refred}. Write
    \begin{eqnarray}\label{eq:case3}
      |S| = \tau_P + \tau_N + \left| P \cap \left[ 1,\, \frac{n}{2} \right] \right| + \left| N \cap \left[ 1,\, \frac{n}{2} \right] \right|.
    \end{eqnarray}
    Since $a \geq \frac{n}{8}$, we have $n-4a \leq \frac{n}{2}$. Thus,
    \begin{eqnarray*}
      \left| N \cap \left[ 1,\, \frac{n}{2} \right] \right| \leq \left| N \cap \left[ 1,\, n-4a \right] \right| + \frac{n}{2} - (n-4a).
    \end{eqnarray*}
    Applying Lemma \ref{lem:4refred} for $a$ and similarly for $b$, and
    inserting into (\ref{eq:case3}), we get
    $|S| \leq 4(a+b) + 1 \leq \left( \frac{33}{32} - 4c \right) n + 1$, which is
    enough for (E1). Moreover, note that no extremal set can arise from Case 3,
    for $n \gg 0$, provided we choose a strictly positive $c$. 
    \\
    \par It now remains to establish (\ref{eq:E1bound}) under the assumptions
    that $a, \, b \leq \frac{n}{6}$ and
    $a+b \geq \left( \frac{33}{128} - c\right) n$,
    for any fixed $0 < c < \frac{1}{128}$. While the choice of $c$ is
    arbitrary, a convenient choice is $c = \frac{3}{512}$. Thus we define
    \\
    \\
    {\sc Case 4:} $\max \{a,\,b\} \leq \frac{n}{6}$ and $a+b \geq \frac{129n}{512}$.
    \\
    \par We will prove that, in Case 4,
    \begin{eqnarray}\label{eq:C4bound}
      |S| \leq \left( \frac{33}{32} - \frac{1}{2048} \right) n + o(n),
    \end{eqnarray}
    which in particular means that no extremal sets can arise in this case
    either, for $n \gg 0$. The rest of Section \ref{sect:kfour} is devoted to
    proving this claim.
    It will require a deeper analysis of trade-offs between $M$ and $T$, and to
    a lesser extent between $Q$ and $T$, in a
    $4$-sum-free set. As preparation for a sequence of seven further lemmas, we
    first need to introduce some more notation.
    \\
    \par If $X=P$, let $O := N$ and{\footnote{$O$ stands for ``other''.}}
    $(u,\,v) := (a,\,b)$. If $X = N$, let $O := P$ and $(u,\,v) = (b,\,a)$.
    Since both $a$ and $b$ are at most $n/6$, by (\ref{eq:M4trade}) we have
    $\kappa(u) = \frac{n/2 - 3u}{4}$ and $\kappa(v) = \frac{n/2 - 3v}{4}$. Let
    \begin{eqnarray*}
      \theta := \frac{n/2 + u}{4}, \;\;\; s_2 := \frac{5u+v}{4}.
    \end{eqnarray*}
    One readily verifies the following:
    \begin{itemize}
    \item
      $\theta = u + \kappa(u)$, for any $u \leq \frac{n}{6}$,
    \item
      Since $a, \, b \geq \frac{n}{8}$ and $a+b > \frac{n}{4}$,
      $tu+v > \frac{t+1}{8} n$, for any $t \geq 1$,
    \item
      $s_2 - \theta > \frac{\kappa(u)}{4}$,
    \item
      $s_2 \leq \frac{n}{4}$, since $a, \, b \leq \frac{n}{6}$.
    \end{itemize}
    In particular, $u \leq \theta < s_2 \leq \frac{n}{4}$, so
    $\mathcal{X} := X \cap M$ can be split into three parts:
    \begin{eqnarray*}
      \mathcal{X}_1 = \mathcal{X} \cap [u,\,\theta], \;\;\;
      \mathcal{X}_2 = \mathcal{X} \cap (\theta, \, s_2], \;\;\;
      \mathcal{X}_3 = \mathcal{X} \cap \left( s_2, \, \frac{n}{4} \right].
    \end{eqnarray*}
    Let $\mu_i = \mu_{i,\,X} = |\mathcal{X}_i|$, so $m_X = \sum_{i=1}^{3} \mu_{i,\,X}$.

    We will also need to split the interval $X \cap Q$ into two parts, namely
    $X \cap Q \cap \left( \frac{n}{4}, \, n-4v\right]$ and
    $X \cap Q \cap \left(n-4v, \, \frac{n}{2} \right]$. Note that this is
    well-defined since $\frac{n}{8} \leq u,\,v \leq \frac{n}{6}$. The sizes of
    the respective parts will be denoted $q_{1} = q_{1,\,X}$ and $\rho=\rho_X$.

    Finally, we let $c_X \in \{1,\,2,\,3,\,4\}$ be the number of residue
    classes mod $4$ met by $X \cap M$. The assumptions of Case 4 disallow
    $c_X = 0$.
    \\
    \par The following four inequalities are all trivial consequences of the
    above definitions:
    \begin{eqnarray}
      \mu_{1,\,X} \leq \frac{c_X \kappa(u)}{4} + O(1), \label{eq:mu1} \\
      \mu_{2,\,X} \leq \frac{c_X (s_2 - \theta)}{4} + O(1), \label{eq:mu2} \\
      \mu_{3,\,X} \leq \frac{c_X (n/4 - s_2)}{4} + O(1), \;\; {\hbox{and $\mu_{3,\,X} = 0$ if $\mathcal{X}_3 = \phi$}}. \label{eq:mu4} \\
      \rho_X \leq 4v - \frac{n}{2}. \label{eq:rho}
    \end{eqnarray}
    Also, in the above notation, (\ref{eq:q4}) becomes
    \begin{eqnarray}\label{eq:Q4qr}
      q_{1,\,X} + \rho_X \leq \frac{g_X}{2} + O(1).
    \end{eqnarray}

    \subsection*{Seven further inequalities.} The above completes the
    preparations for the deeper analysis of Case 4, and we proceed now with the
    proofs of seven additional lemmas. In what follows, when a subscript $_X$
    refers to $X \in \{P,\,N\}$, the subscript $_O$ refers to the other member
    of $\{P, \, N\}$, thus to the ``other side'' of a $4$-sum-free set $S$.  

    \begin{lemma}\label{lem:R}
      For $X \in\{P, \, N\}$,
      \begin{eqnarray}\label{eq:R}
        g_X \geq l_O + m_O + q_{1,\,O}, \;\;\; {\hbox{where $q_{1,\,O} = \left| O \cap \left( \frac{n}{4}, \, n-4u \right] \right|$}}.
      \end{eqnarray}
    \end{lemma}

    \begin{proof}
      Applying Lemma \ref{lem:4refred} at $x=u \in X$, we get
      \begin{eqnarray*}
        \tau_X + |O \cap [1,\,n-4u]| \leq |T| \;\; \Leftrightarrow \;\; g_X \geq |O \cap [1,\,n-4u]|.
      \end{eqnarray*}
      But
      \begin{eqnarray*}
        |O \cap [1,\,n-4u]| = |O \cap L| + |O \cap M| + \left| O \cap \left( \frac{n}{4}, \, n-4u \right] \right| = l_O + m_O + q_{1,\,O}.
        \end{eqnarray*}
    \end{proof}
    
    \begin{lemma}\label{lem:DO}
      For $X \in \{P,\,N\}$,
      \begin{eqnarray}\label{eq:DO}
        g_X \geq c_O \mu_{1,\,X} + m_O - O(1).
      \end{eqnarray}
\end{lemma}

\begin{proof}
  Consider the set $A = 4\mathcal{X}_1 + (O \cap M)$. It avoids $X$, by parts
  (iii) and (iv) of Lemma \ref{lem:elem}. It is a subset of
  $\left[ 4u+v,\,4\theta + \frac{n}{4} \right] = \left[4u+v, \, u + \frac{3n}{4} \right]$ and hence of $T$, so $g_X \geq |A|$. For each residue class
  $r \; ({\hbox{mod $4$}})$ met by $O \cap M$, using (\ref{eq:Intsumset}) we
  have
  \begin{eqnarray*}
    |A_r| = |4\mathcal{X}_1 + (O \cap M)_r| \geq |\mathcal{X}_1| + |(O \cap M)_r| - 1 = \mu_{1,\,X} + |(O \cap M)_r| - 1.
  \end{eqnarray*}
  Thus,
  \begin{eqnarray*}
    g_X \geq |A| = \sum_r |A_r| = c_O (\mu_{1,\,X} - 1) + |O \cap M| \geq c_O \mu_{1,\,X} + m_O - 4.
  \end{eqnarray*}
\end{proof}

\begin{lemma}\label{lem:DS}
  For $X \in \{P,\,N\}$,
  \begin{eqnarray}\label{eq:DS}
    g_X \geq (4+c_O) \mu_{1,\,X} + 4\mu_{2,\,X} + m_O - c_X \kappa(u) - o(n).
    \end{eqnarray}
  \end{lemma}

  \begin{proof}
    Let $Z = \mathcal{X}_1 \cup \mathcal{X}_2 \subseteq [u,\,s_2]$ and
    $Z^{*} = Z-u$. The set $4Z - Z$ avoids $X$ and
    $\max (4Z - Z) \leq 4s_2 - u = 4u+v$. Hence, $|(4Z - Z) \cap A| \leq 1$,
    where $A$ is the set in Lemma \ref{lem:DO}, and so
    \begin{eqnarray}\label{eq:DS1}
      g_X \geq c_O \mu_{1,\,X} + m_O - O(1) + g', \;\;\; {\hbox{where $g' = \left| (4Z-Z) \cap \left( \frac{n}{2}, \, 4u+v \right]\right|$}}.
    \end{eqnarray}
    The numbers in $4Z - Z$ above $n/2$ are in one-to-one correspondence with
    those of $Z^{*}$ above $K = n/2  - 3u = 4\kappa(u)$, where $Z^{*} = Z - u$.
    Hence, by Corollary \ref{cor:dilates}(ii),
    \begin{eqnarray*}
      g' = P(Z^{*}, \, K) \geq 4|Z^{*}| - c_X \kappa(u) - o(n) = 4(\mu_{1,\,X} + \mu_{2,\,X}) - c_X \kappa(u) - o(n).
    \end{eqnarray*}
    Inserting this into (\ref{eq:DS1}) yields (\ref{eq:DS}).
  \end{proof}

  \begin{remark}\label{rem:DS}
    Though the proof of (\ref{eq:DS}) might suggest it is a stronger inequality
    than (\ref{eq:DO}), this is not necessarily the case since the lower bound
    obtained from Corollary \ref{cor:dilates}(ii) could, in principle, be
    negative. Similar remarks apply to the comparison between (\ref{eq:DO}) and
    (\ref{eq:DW}) below.
    \end{remark}
  
\begin{lemma}\label{lem:DW}
  For $X \in \{P, \, N\}$,
  \begin{eqnarray}\label{eq:DW}
    g_X \geq (c_O + \beta) \mu_{1,\,X} + 4\mu_{2,\,X} + m_O - \kappa(u) - o(n), \\ {\hbox{where $\beta = \left\{ \begin{array}{lr} 1, & {\hbox{if $\mathcal{X}_2 \neq \phi$ and $\max \mathcal{X}_2 \geq \theta + \frac{\kappa (u)}{4}$}}, \\ 0, & {\hbox{otherwise}}. \end{array} \right.$}} \nonumber
  \end{eqnarray}
  \end{lemma}

\begin{proof}
  If $\mathcal{X}_2$ is empty, then $\mu_{2,\,X} = 0$ and (\ref{eq:DW}) follows
  from (\ref{eq:DO}). Otherwise, still using the proof of Lemma \ref{lem:DO},
  it suffices to find at least $\beta \cdot \mu_{1,\,X} + 4 \mu_{2,\,X} - \kappa(u) - o(n)$ elements of $X^c$ in the interval
  $\left( \frac{n}{2}, \, 4u+v \right]$.

  Let $m_2 = \min \mathcal{X}_2$ and $\mathcal{X}^{*}_{2} = \mathcal{X}_2 - m_2$.
  The set $4\mathcal{X}_2 - \mathcal{X}_2$ is a subset of the set $4Z-Z$ from
  Lemma \ref{lem:DS}, so it is disjoint from $X$ and its elements don't exceed
  $4u+v$. Those exceeding $n/2$ are in one-to-one correspondence with those of
  $4\mathcal{X}^{*}_{2} - \mathcal{X}^{*}_2$ exceeding
  $n/2 - 3m_2 < n/2 - 3\theta = \kappa(u)$. By Corollary \ref{cor:dilates}(i),
  there are at least $4\mu_{2,\,X} - \kappa(u) - o(n)$ such elements. This
  already suffices when $\beta = 0$, whereas if $\beta = 1$ we still need to
  find a further $\mu_{1,\,X} - o(n)$ elements of $X^c$ in the interval
  $\left( \frac{n}{2}, \, 4u+v \right]$.

  Let $M_2 = \max \mathcal{X}_2$. If $\beta = 1$, then
  $M_2 \geq \theta + \frac{\kappa (u)}{4}$. The set $4M_2 - \mathcal{X}_1$ has
  size $\mu_{1,\,X}$, is disjoint from $X$ and obviously disjoint from
  $4\mathcal{X}_2 - \mathcal{X}_2$. It is also a subset of $4Z - Z$, so none of
  its elements exceed $4u+v$. But if $\beta = 1$, then all of its elements are
  at least $4\left( \theta + \frac{\kappa (u)}{4} \right) - \theta = \kappa (u) + 3\theta = \frac{n}{2}$. So at most one element in this set isn't above $n/2$,
  and we are done. 
\end{proof}

The next lemma is a strengthening of Lemma \ref{lem:M4trade}, up to an $o(n)$
term.

\begin{lemma}\label{lem:WS}
  For $X \in \{P,\,N\}$,
  \begin{eqnarray}\label{eq:WS}
    g_X \geq 4 m_X + \delta \cdot l_X - c_X \kappa (u) - o(n),
    \\
    {\hbox{where $\delta = \left\{ \begin{array}{lr} 1, & {\hbox{if $\mathcal{X}_2 \cup \mathcal{X}_3 \neq \phi$}}, \\ 0, & {\hbox{otherwise}}. \end{array} \right.$}} \nonumber
  \end{eqnarray}
\end{lemma}

\begin{proof}
  Let $A^{*} = \mathcal{X} - u$ be as in the proof of Lemma \ref{lem:M4trade},
  and repeat the proof of that lemma, but applying part (ii) of Corollary
  \ref{cor:dilates} instead of part (i). This already gives (\ref{eq:WS}) when
  $\delta = 0$.

  If $\delta = 1$, it suffices to find an additional $l_X$ elements of
  $X^c \cap T$, disjoint from $4\mathcal{X} - \mathcal{X}$. Let $m^{*}$ be the
  largest element in $\mathcal{X} = X \cap M$. The set $4m^{*} - (X \cap L)$
  has size $l_X$, and is evidently disjoint from both $X$ and
$4\mathcal{X} - \mathcal{X}$. If $\delta = 1$, then $m^{*} > \theta$ and every
element of $4m^{*} - (X \cap L)$ is strictly greater than
$4\theta - \frac{n}{8} = \frac{n}{2} + \left(u - \frac{n}{8} \right) \geq \frac{n}{2}$. So this set lies entirely in $X^c \cap T$, as desired.
\end{proof}

\begin{lemma}\label{lem:WW}
  Let $X \in \{P,\,N\}$ and suppose $\mathcal{X}_2 \cup \mathcal{X}_3 \neq \phi$. Then
  \begin{eqnarray}\label{eq:WW}
    g_X \geq 4(\mu_{2,\,X} + \mu_{3,\,X}) + l_X + \gamma \cdot \mu_{1,\,X} - \kappa(u) - o(n), \\
    {\hbox{where $\gamma = \left\{ \begin{array}{lr} 1, & {\hbox{if $\max (\mathcal{X}_2 \cup \mathcal{X}_3) \geq \theta + \frac{\kappa (u)}{4}$}}, \\ 0, & {\hbox{otherwise}}. \end{array} \right.$}} \nonumber
  \end{eqnarray}
\end{lemma}

\begin{proof}
  Let $m_{*}$ and $m^{*}$ be the smallest and largest elements, respectively, of
  $\mathcal{X}_2 \cup \mathcal{X}_3$. Let
  $Z = \mathcal{X}_2 \cup \mathcal{X}_3$. Applying Corollary
  \ref{cor:dilates}(i) to $Z-m_{*}$ and noting that
  $\frac{n}{2} - 3m_{*} < \frac{n}{2} - 3\theta = \kappa(u)$, we find that at
  least $4(\mu_{2,\,X} + \mu_{3,\,X}) - \kappa(u) - o(n)$ elements of $4Z - Z$
  are in $X^c \cap T$.

  Then, as in the proof of Lemma \ref{lem:WS}, one easily checks that:
  \begin{itemize}
  \item
    If $\gamma = 0$, $4m^{*} - x$ is greater than every element in $4Z - Z$, for
    $x \in X \cap L$.
  \item
    If $\gamma = 1$, $4m^{*} - x$ is greater than every element in $4Z - Z$, for
    $x \in X \cap (L \cup \mathcal{X}_1)$.
  \end{itemize}
  Adding these elements to those of $4Z - Z$ yields (\ref{eq:WW}).
  \end{proof}

  \begin{lemma}\label{lem:U1}
    Let $X \in \{P,\,N\}$ and suppose $\mathcal{X}_3 \neq \phi$. Then
    \begin{eqnarray}\label{eq:U1}
      g_X \geq 5\mu_{3,\,X} + \mu_{1,\,X} + \mu_{2,\,X} + l_X - O(1).
      \end{eqnarray}
\end{lemma}

\begin{proof}
  Every element of $4\mathcal{X}_3 - \mathcal{X}_3$ strictly exceeds
  $4s_2 - \frac{n}{4} = (5u+v) - \frac{n}{4} > \frac{5+1}{8} n - \frac{n}{4} = \frac{n}{2}$. Thus all of $4\mathcal{X}_3 - \mathcal{X}_3$ lies in $X^c \cap T$,
  and this set has at least $5\mu_{3,\,X} - O(1)$ elements, by Shakan's result
  \cite{S}. Then, every element of
  $4m^{*} - ((X \cap L) \cup \mathcal{X}_1 \cup \mathcal{X}_2)$, where
  $m^{*} = \max \mathcal{X}_3$, strictly exceeds every element of
  $4\mathcal{X}_3 - \mathcal{X}_3$, which gives an additional
  $l_X + \mu_{1,\,X} + \mu_{2,\,X}$ elements of $X^c \cap T$.
\end{proof}

The following table summarizes the inequalities that we have proven above, and
which will be combined in various ways to complete the analysis of Case 4. In
each inequality where $g_X$ appears, we write it instead in terms of
$\tau_X = |T| - g_X$.

\begin{table}[ht!]
  \begin{center}
    \begin{tabular}{|c|c|} \hline
      (\ref{eq:l4a}) & $l_P + l_N + \frac{1}{4} (\tau_P + \tau_N) \leq \frac{n}{4} + O(1)$. \\ \hline
      (\ref{eq:mu1}) & $\mu_{1,\,X} \leq \frac{c_X \kappa (u)}{4} + O(1)$. \\ \hline
       (\ref{eq:mu2}) & $\mu_{2,\,X} \leq \frac{c_X (s_2 - \theta)}{4} + O(1)$. \\ \hline
       (\ref{eq:mu4}) & $\mu_{3,\,X} \leq \frac{c_X (n/4 - s_2)}{4} + O(1); \;\; \mu_{3,\,X} = 0, \, {\hbox{if $\mathcal{X}_3 = \phi$}}$. \\ \hline
        (\ref{eq:rho}) & $\rho_X \leq 4v - \frac{n}{2}$. \\ \hline
        (\ref{eq:Q4qr}) & $q_{1,\,X} + \rho_X + \frac{\tau_X}{2} \leq \frac{n}{4} + O(1)$. \\ \hline
         (\ref{eq:R}) & $l_O + m_O + q_{1,\,O} + \tau_X \leq \frac{n}{2} + 1$. \\ \hline
         (\ref{eq:DO}) & $c_O \mu_{1,\,X} + m_O + \tau_X \leq \frac{n}{2} + O(1)$. \\ \hline
          (\ref{eq:DS}) & $(4 + c_O)\mu_{1,\,X} + 4\mu_{2,\,X} + m_O + \tau_X \leq \frac{n}{2} + c_X \kappa(u) + o(n)$. \\ \hline
          (\ref{eq:DW}) & $(c_O + \beta) \mu_{1,\,X} + 4 \mu_{2,\,X} + m_O + \tau_X \leq \frac{n}{2} + \kappa(u) + o(n)$. \\ \hline
           (\ref{eq:WS}) & $\delta \cdot l_X + 4m_X + \tau_X \leq \frac{n}{2} + c_X \kappa(u) + o(n)$. \\ \hline
           (\ref{eq:WW}) & $l_X + \gamma \cdot \mu_{1,\,X} + 4(\mu_{2,\,X} + \mu_{3,\,X}) + \tau_X \leq \frac{n}{2} + \kappa(u) + o(n), \; {\hbox{if $\mathcal{X}_2 \cup \mathcal{X}_3 \neq \phi$}}$. \\ \hline
            (\ref{eq:U1}) & $l_X + \mu_{1,\,X} + \mu_{2,\,X} + 5\mu_{3,\,X} + \tau_X \leq \frac{n}{2} + O(1), \; {\hbox{if $\mathcal{X}_3 \neq \phi$}}$. \\ \hline
          \end{tabular}
        \end{center}
  \caption{List of the inequalities which will be combined to complete the analysis of Case 4.}
  \label{tab:one}
  \end{table}
  
  The size of a $4$-sum-free set $S \subseteq [-n,\,n]$ can be expressed as a
  sum of 14 quantities:
  \begin{eqnarray}\label{eq:Sdecomp}
    |S| = \sum_{X \in \{P,\,N\}} l_X + (\mu_{1,\,X} + \mu_{2,\,X} + \mu_{3,\,X}) + (q_{1,\,X} + \rho_X) + \tau_X.
  \end{eqnarray}
  Note that the left-hand side of each inequality in Table \ref{tab:one} is
  some affine combination of these 14 quantities. Since we're in Case 4, the
  point $\left(\frac{a}{n}, \, \frac{b}{n} \right)$ belongs to the planar region
  \begin{eqnarray*}
    \mathcal{R} = \left\{ (x,\,y) : x, \, y \in \left[ \frac{1}{8}, \, \frac{1}{6} \right], \; x+y \geq \frac{129}{512} \right\},
  \end{eqnarray*}
  which is a pentagon with vertices at, reading anti-clockwise,
  \begin{eqnarray*}
    \left( \frac{1}{8}, \, \frac{65}{512} \right), \;
    \left( \frac{65}{512}, \, \frac{1}{8} \right), \;
    \left( \frac{1}{6}, \, \frac{1}{8} \right), \;
    \left( \frac{1}{6}, \, \frac{1}{6} \right), \;
    \left( \frac{1}{8}, \, \frac{1}{6} \right).
  \end{eqnarray*}
  A \emph{certificate} is a list of non-negative weights on the 13 inequalities
  in Table \ref{tab:one} such that, in the weighted sum, the coefficient of
  each of the 14 quantities in (\ref{eq:Sdecomp}) is at least $1$. Then $|S|$
  will be at most the weighted sum of the right-hand sides of the 13
  inequalities. The latter weighted sum will always be an affine function of
  $a$ and $b$, plus an $o(n)$ term. Hence, to complete Case 4, it will suffice
  to show that this affine function is at most
  $\left( \frac{33}{32} - \frac{1}{2048} \right) n$ at each of the 5 vertices
  of the region $\mathcal{R}$. 

  We cannot, however, use a single certificate, because the 13 inequalities
  also depend on the values of $c_P, \, c_N \in \{1,\,2,\,3,\,4\}$ plus whether
  or not each of the sets $\mathcal{X}_2, \, \mathcal{X}_3$, $X \in \{P,\,N\}$
  is empty. This gives a total of $4 \cdot 4 \cdot 2^4 = 256$ possible
  ``states''. It turns out that 14 different certificates will suffice to cover
  all the states. The first of these (S1) is presented in full below, to show
  for the reader how things work. Six more (S2)-(S7) are presented in the
  Appendix. The remaining 7 certificates are the ``mirror images'' of these
  seven, got by interchanging $P$ and $N$. We will leave it to the reader to
  verify that the definitions of the 14 certificates cover all 256 states
  (with some overlap).

  Before presenting (S1), we note a number of further relevant points:
  \begin{itemize}
  \item
    There are actually two versions of each of the 13 inequalities in Table
    \ref{tab:one}, one for $X=P$ and one for $X=N$.
  \item
    For $X \in \{P,\,N\}$, $c_X$ only ever appears on the right-hand side of an
    inequality in the table, whereas $c_O$ only ever appears on the left-hand
    side. This means that any inequality involving either $c_X$ or $c_O$ will
    remain valid if we replace $c_X$ by any upper bound for it, and $c_O$ by
    any lower bound for it.
  \item
    None of the 14 certificates take account of whether or not $\mathcal{X}_2$
    is empty, $X \in \{P,\,N\}$, so there are effectively 64 states rather than
    256, divided into 14 overlapping groups.
  \item
    Inequality (\ref{eq:DW}) is only ever applied with $\beta = 0$.
  \end{itemize}
  


\subsection*{Certificate (S1).} This certificate applies when
$c_P = 1$, $1 \leq c_N \leq 4$, $\mathcal{P}_3 \neq \phi$, $\mathcal{N}_3 \neq \phi$. On the right-hand sides we use the upper ends $c_P = 1$, $c_N = 4$; as
the coefficient $c_O$ in (\ref{eq:DO}), (\ref{eq:DS}), (\ref{eq:DW}) we use
the lower ends, i.e.~$c_N = 1$ for $X = P$ and $c_P = 1$ for $X = N$.
\par Each row of the table below is one of the inequalities of Table
\ref{tab:one}, for the stated $X \in \{P,\,N\}$, multiplied by its weight; the
columns give the coefficients of the 14 quantities of (\ref{eq:Sdecomp}) on
the left-hand side (before multiplication by the weight) and the right-hand
side, written as an affine function of $n$, $a$ and $b$ (recall that
$(u,\,v) = (a,\,b)$ for $X=P$ and $(u,\,v) = (b,\,a)$ for $X=N$), with the
$O(1)$ and $o(n)$ terms omitted. Here $\mathcal{P}_3, \, \mathcal{N}_3$ denote
$\mathcal{X}_3$ for $X = P, \, N$. Where $\mathcal{X}_3 \neq \phi$, we use
(\ref{eq:WS}) with $\delta = 1$ and (\ref{eq:WW}) with $\gamma = 1$: indeed
$\max (\mathcal{X}_2 \cup \mathcal{X}_3) > s_2 > \theta + \frac{\kappa(u)}{4}$.
Where $\mathcal{X}_3$ may be empty, we use (\ref{eq:WS}) with $\delta = 0$, and
neither (\ref{eq:WW}) nor (\ref{eq:U1}). The last row is the weighted sum:
every coefficient in it is at least $1$, so $|S|$ is at most its right-hand
side, up to $o(n)$.
\begin{center}
{\scriptsize \setlength{\tabcolsep}{2pt} \renewcommand{\arraystretch}{1.4}
\resizebox{\textwidth}{!}{%
\begin{tabular}{|c|c|c|ccccccc|ccccccc|c|} \hline
Ineq. & $X$ & Weight & $l_P$ & $\mu_{1,\,P}$ & $\mu_{2,\,P}$ & $\mu_{3,\,P}$ & $q_{1,\,P}$ & $\rho_P$ & $\tau_P$ & $l_N$ & $\mu_{1,\,N}$ & $\mu_{2,\,N}$ & $\mu_{3,\,N}$ & $q_{1,\,N}$ & $\rho_N$ & $\tau_N$ & Right-hand side \\ \hline
(\ref{eq:l4a}) &  & $\frac{3}{4}$ & $1$ &  &  &  &  &  & $\frac{1}{4}$ & $1$ &  &  &  &  &  & $\frac{1}{4}$ & $\frac{n}{4}$ \\
(\ref{eq:Q4qr}) & $P$ & $1$ &  &  &  &  & $1$ & $1$ & $\frac{1}{2}$ &  &  &  &  &  &  &  & $\frac{n}{4}$ \\
(\ref{eq:R}) & $P$ & $\frac{1}{16}$ &  &  &  &  &  &  & $1$ & $1$ & $1$ & $1$ & $1$ & $1$ &  &  & $\frac{n}{2}$ \\
(\ref{eq:WS}) & $P$ & $\frac{11}{48}$ & $1$ & $4$ & $4$ & $4$ &  &  & $1$ &  &  &  &  &  &  &  & $\frac{5n}{8} - \frac{3a}{4}$ \\
(\ref{eq:U1}) & $P$ & $\frac{1}{48}$ & $1$ & $1$ & $1$ & $5$ &  &  & $1$ &  &  &  &  &  &  &  & $\frac{n}{2}$ \\
(\ref{eq:rho}) & $N$ & $\frac{1}{24}$ &  &  &  &  &  &  &  &  &  &  &  &  & $1$ &  & $4a - \frac{n}{2}$ \\
(\ref{eq:Q4qr}) & $N$ & $\frac{23}{24}$ &  &  &  &  &  &  &  &  &  &  &  & $1$ & $1$ & $\frac{1}{2}$ & $\frac{n}{4}$ \\
(\ref{eq:DS}) & $N$ & $\frac{5}{48}$ &  & $1$ & $1$ & $1$ &  &  &  &  & $5$ & $4$ &  &  &  & $1$ & $n - 3b$ \\
(\ref{eq:WS}) & $N$ & $\frac{1}{16}$ &  &  &  &  &  &  &  & $1$ & $4$ & $4$ & $4$ &  &  & $1$ & $n - 3b$ \\
(\ref{eq:WW}) & $N$ & $\frac{1}{12}$ &  &  &  &  &  &  &  & $1$ & $1$ & $4$ & $4$ &  &  & $1$ & $\frac{5n}{8} - \frac{3b}{4}$ \\
(\ref{eq:U1}) & $N$ & $\frac{1}{12}$ &  &  &  &  &  &  &  & $1$ & $1$ & $1$ & $5$ &  &  & $1$ & $\frac{n}{2}$ \\
\hline
Sum &  &  & $1$ & $\frac{25}{24}$ & $\frac{25}{24}$ & $\frac{9}{8}$ & $1$ & $1$ & $1$ & $\frac{25}{24}$ & $1$ & $\frac{55}{48}$ & $\frac{17}{16}$ & $\frac{49}{48}$ & $1$ & $1$ & $\frac{141n}{128} - \frac{a}{192} - \frac{9b}{16}$ \\ \hline
\end{tabular}}}
\end{center}
At the five vertices of $\mathcal{R}$, listed as above, the right-hand side of
the last row equals $\frac{25301n}{24576}$, $\frac{101311n}{98304}$,
$\frac{1187n}{1152}$, $\frac{145n}{144}$, $\frac{1547n}{1536}$. The maximum is
$\frac{101311n}{98304} < \left( \frac{33}{32} - \frac{1}{2048} \right) n$.

\vspace{3mm}

\section{{\bf Case of $k=4$: Full proof of Theorem \ref{thm:B}}}\label{sect:stability}

We use the same notation here as in Section \ref{sect:kfour}, introducing
further notation as necessary.

The proof of Theorem \ref{thm:B} will be completed in three steps.

The first step is to reduce to Cases 1--4 from the previous section. There we
first showed that a $4$-sum-free subset of $[-n,\,n]$ must satisfy either (E1)
or (E2). If (E2) holds then, since we have now established (\ref{eq:Off}),
\begin{eqnarray*}
  |S| \leq \frac{n}{2} + O(1) + \frac{33}{32}\left(\frac{n}{2}\right) + o(n) = \left( \frac{33}{32} - \frac{1}{64} \right) n + o(n),
\end{eqnarray*}
which means $S$ cannot have maximum size, when $n \gg 0$. The same applies when
either side of $S$ is exceptional, by (\ref{eq:excbd}). This shows that an
extremal $S$ must fall under Cases 1--4, for $n \gg 0$.
\\
\par The second step is to eliminate Cases 2--4, as follows. 

\begin{theorem}\label{thm:cases2to4}
  There is an absolute constant $c > 0$ such that, if $S$ is a $4$-sum-free
  subset of $[-n,\,n]$ for which $P \cap M$ and $N \cap M$ are both non-empty,
  then $|S| \leq \left( \frac{33}{32} - c \right) n + o(n)$.
\end{theorem}

\begin{proof}
  The previous discussion proves the existence of such a $c > 0$, unless we are
  in one of Cases 1--4 from Section \ref{sect:kfour}. If both $P \cap M$ and
  $N \cap M$ are non-empty, then Case 1 is excluded. The existence of an
  absolute $c > 0$ in Cases 3 and 4 was already established in Section
  \ref{sect:kfour}. In Case 2, we used the trivial bound
  $\kappa(u) \leq \frac{n}{32}$, which only gives $c=0$. So it remains to
  conduct a more careful analysis of Case 2. Note that, if both $a$ and $b$ are
  greater than $\frac{n}{6}$, then $\kappa(a) = \kappa(b) = 0$ and
  $|S| \leq n + o(n)$, so we may assume exactly one of $a$ and $b$ is greater
  than $\frac{n}{6}$. Without loss of generality, this is $b = \min N \cap M$. 

  The strategy now is to go through the list of inequalities in Table
  \ref{tab:one} and observe that many of them still hold in the current
  setting, enough so that we can still construct certificates similar to those
  in Section \ref{sect:kfour} which cover all possibilities. In fact, it turns
  out that the certificates are much simpler now. But we begin by going through
  Table \ref{tab:one}.
  \begin{itemize}
  \item
    Eq. (\ref{eq:l4a}) makes no reference to the values of $a$ and $b$, and
    still holds.
  \item
    Eq. (\ref{eq:mu1}) is a simple count and still holds. Note that now
    $\kappa(u) = 0$ when $X = N$. We will only apply it for $X=P$.
    Eq. (\ref{eq:mu2}) is also a simple count and still holds, but we will not
    need it. Eq. (\ref{eq:mu4}) may no longer hold since $s_2 > \frac{n}{4}$ is
    possible when $v > \frac{n}{6}$. We will not need it either.
  \item
    The decomposition of $Q$ into two parts relied on having
    $n-4v> \frac{n}{4}$, and thus doesn't make sense if $v > \frac{3n}{16}$. So
    we do not use (\ref{eq:rho}) and (\ref{eq:Q4qr}). However, the original
    (\ref{eq:q4}) still holds, and we will use that.
  \item
    We will not use (\ref{eq:R}).
  \item
    The proof of Lemma \ref{lem:DO} goes through as long as
    $u \leq \frac{n}{6}$: note, in particular, that
    $4u + v \leq 4\left( \frac{n}{6} \right) + \frac{n}{4} \leq \frac{11n}{12} < n$. So we can still use (\ref{eq:DO}) for $X = P$. The same applies to
    (\ref{eq:DS}) and to (\ref{eq:DW}) for $\beta = 0$, which was also the only
    form used in Section \ref{sect:kfour}. 
  \item
    We want to apply (\ref{eq:WS}) for $X=N$ and $\delta = 1$. The first part
    of the proof of Lemma \ref{lem:WS} only appealed to that of Lemma
    \ref{lem:M4trade} and thus goes through for any value of $u$. For the case
    $\delta = 1$ to make sense, we require
    $\mathcal{X}_2 \cup \mathcal{X}_3 \neq \phi$. This is now automatic when
    $X = N$, because then $u = b > \frac{n}{6}$, $\mathcal{X}_1$ is empty and
    the point $b = u$ is itself in $\mathcal{X}_2$. The second part of the
    proof then just required $u \geq \frac{n}{8}$, which automatically holds.
    Hence, Lemma \ref{lem:WS} still holds for $X=N$, and in the current
    setting, $\kappa(u) = \kappa(b) = 0$. We will thus apply (\ref{eq:WS}) in
    the form $g_N \geq 4m_N + l_N - o(n)$.
  \item
    Finally, we will only use (\ref{eq:WW}) and (\ref{eq:U1}) under the
    assumption that $\mathcal{X}_3$ is non-empty. In that case, the proofs of
    both go through provided $u \leq \frac{n}{6}$ and thus still apply when
    $X = P$. We can take $\gamma = 1$ in (\ref{eq:WW}).
  \end{itemize}

  The inequalities which we will use are summarized in the following table. 

  \begin{table}[ht!]
  \begin{center}
    \begin{tabular}{|c|c|} \hline
      (\ref{eq:l4a}) & $l_P + l_N + \frac{1}{4} (\tau_P + \tau_N) \leq \frac{n}{4} + O(1)$. \\ \hline
      (\ref{eq:mu1})$_P$ & $\mu_{1,\,P} \leq \frac{c_P \kappa (u)}{4} + O(1)$. \\ \hline
        (\ref{eq:q4}) & $q_X + \frac{\tau_X}{2} \leq \frac{n}{4} + O(1), \;\; X \in \{P,\,N\}$. \\ \hline
         (\ref{eq:DO})$_P$ & $c_N \mu_{1,\,P} + m_N + \tau_P \leq \frac{n}{2} + O(1)$. \\ \hline
          (\ref{eq:DS})$_P$ & $(4 + c_N)\mu_{1,\,P} + 4\mu_{2,\,P} + m_N + \tau_P \leq \frac{n}{2} + c_P \kappa(u) + o(n)$. \\ \hline
          (\ref{eq:DW})$_P$ & $c_N  \mu_{1,\,P} + 4 \mu_{2,\,P} + m_N + \tau_P \leq \frac{n}{2} + \kappa(u) + o(n) \;\; (\beta = 0)$. \\ \hline
           (\ref{eq:WS})$_N$ & $l_N + 4m_N + \tau_N \leq \frac{n}{2} + o(n) \;\; (\delta = 1)$. \\ \hline
           (\ref{eq:WW})$_P$ & $l_P + \mu_{1,\,P} + 4(\mu_{2,\,P} + \mu_{3,\,P}) + \tau_P \leq \frac{n}{2} + \kappa(u) + o(n), \; {\hbox{if $\mathcal{P}_3 \neq \phi$}}$. \\ \hline
            (\ref{eq:U1})$_P$ & $l_P + \mu_{1,\,P} + \mu_{2,\,P} + 5\mu_{3,\,P} + \tau_P \leq \frac{n}{2} + O(1), \; {\hbox{if $\mathcal{P}_3 \neq \phi$}}$. \\ \hline
          \end{tabular}
        \end{center}
  \caption{List of the inequalities which will be combined to complete the analysis of Case 2.}
  \label{tab:two}
\end{table}

We only need two certificates, depending on whether or not $\mathcal{P}_3$ is
empty: note that a necessary condition for non-emptiness is that
$s_2 < \frac{n}{4} \Leftrightarrow 5a + b < n$. When $\mathcal{P}_3$ is empty,
instead of (\ref{eq:WW})$_P$ and (\ref{eq:U1})$_P$ we will use the trivial
inequality
\begin{eqnarray}\label{eq:Z}
  \mu_{3,\,P} \leq 0, \;\; {\hbox{if $\mathcal{P}_3 = \phi$}}.
\end{eqnarray}
We apply the inequalities in Table \ref{tab:two} with $c_P \leq 4$ and
$c_N \geq 1$.


\subsection*{Certificate (S1)'.} Applies when $\mathcal{P}_3 = \phi$.

\begin{center}
{\scriptsize \setlength{\tabcolsep}{2pt} \renewcommand{\arraystretch}{1.4}
\resizebox{\textwidth}{!}{%
\begin{tabular}{|c|c|cccccc|cccc|c|} \hline
Inequality & Weight & $l_P$ & $\mu_{1,\,P}$ & $\mu_{2,\,P}$ & $\mu_{3,\,P}$ & $q_{P}$ & $\tau_P$ & $l_N$ & $m_{N}$ &  $q_{N}$ & $\tau_N$ & Right-hand side \\ \hline
(\ref{eq:Z}) & $1$ &  &  &  & $1$ &  &  &  &  &  &  & $0$ \\
(\ref{eq:q4})$_P$ & $1$ &  &  &  &  & $1$ & $\frac{1}{2}$ &  &  &  &  & $\frac{n}{4} + O(1)$ \\
(\ref{eq:q4})$_N$ & $1$ &  &  &  &  &  &  &  &  & $1$ & $\frac{1}{2}$ & $\frac{n}{4} + O(1)$ \\
(\ref{eq:l4a}) & $1$ & $1$ &  &  &  &  & $\frac{1}{4}$ & $1$ &  &  & $\frac{1}{4}$ & $\frac{n}{4} + O(1)$ \\
(\ref{eq:DS})$_P$ & $\frac{3}{16}$ &  & $5$ & $4$ &  &  & $1$ &  & $1$ &  &  & $n - 3a + o(n)$ \\
(\ref{eq:DW})$_P$ & $\frac{1}{16}$ &  & $1$ & $4$ &  &  & $1$ &  & $1$ &  &  & $\frac{5n}{8} - \frac{3a}{4} + o(n)$ \\
(\ref{eq:WS})$_N$ & $\frac{1}{4}$ &  &  &  &  &  &  & $1$ & $4$ &  & $1$ & $\frac{n}{2} + o(n)$ \\
\hline
Weighted sum &  & $1$ & $1$ & $1$ & $1$ & $1$ & $1$ & $\frac{5}{4}$ & $\frac{5}{4}$ & $1$ & $1$ & $\frac{141}{128} n - \frac{39}{64} a + o(n)$ \\ \hline
\end{tabular}}}
\end{center}
Since $a \geq \frac{n}{8}$, the weighted sum on the right-hand side is at most $\left( \frac{33}{32} - c \right) n + o(n)$, for $c = \frac{3}{512}$.

\newpage

\subsection*{Certificate (S2)'.} Applies when $\mathcal{P}_3 \neq \phi$.

\begin{center}
{\scriptsize \setlength{\tabcolsep}{2pt} \renewcommand{\arraystretch}{1.4}
\resizebox{\textwidth}{!}{%
\begin{tabular}{|c|c|cccccc|cccc|c|} \hline
Inequality & Weight & $l_P$ & $\mu_{1,\,P}$ & $\mu_{2,\,P}$ & $\mu_{3,\,P}$ & $q_{P}$ & $\tau_P$ & $l_N$ & $m_{N}$ &  $q_{N}$ & $\tau_N$ & Right-hand side \\ \hline
(\ref{eq:q4})$_P$ & $1$ &  &  &  &  & $1$ & $\frac{1}{2}$ &  &  &  &  & $\frac{n}{4} + O(1)$ \\
(\ref{eq:q4})$_N$ & $1$ &  &  &  &  &  &  &  &  & $1$ & $\frac{1}{2}$ & $\frac{n}{4} + O(1)$ \\
(\ref{eq:l4a}) & $\frac{3}{4}$ & $1$ &  &  &  &  & $\frac{1}{4}$ & $1$ &  &  & $\frac{1}{4}$ & $\frac{n}{4} + O(1)$ \\
(\ref{eq:mu1})$_P$ & $\frac{7}{16}$ &  & $1$ &  &  &  &  &  &  &  &  & $\frac{n}{8} - \frac{3a}{4} + O(1)$ \\
(\ref{eq:DS})$_P$ & $\frac{1}{16}$ &  & $5$ & $4$ &  &  & $1$ &  & $1$ &  &  & $n - 3a + o(n)$ \\
(\ref{eq:WS})$_N$ & $\frac{5}{16}$ &  &  &  &  &  &  & $1$ & $4$ &  & $1$ & $\frac{n}{2} + o(n)$ \\
  (\ref{eq:WW})$_P$ & $\frac{3}{16}$ & $1$ & $1$ & $4$ & $4$ &  & $1$ &  &  &  & & $\frac{5n}{8} - \frac{3a}{4} + o(n)$ \\
  (\ref{eq:U1})$_P$ & $\frac{1}{16}$ & $1$ & $1$ & $1$ & $5$ &  & $1$ &  &  &  &  & $\frac{n}{2} + O(1)$ \\
  \hline
Weighted sum &  & $1$ & $1$ & $\frac{17}{16}$ & $\frac{17}{16}$ & $1$ & $1$ & $\frac{17}{16}$ & $\frac{21}{16}$ & $1$ & $1$ & $\frac{71}{64} n - \frac{21}{32} a + o(n)$ \\ \hline
\end{tabular}}}
\end{center}
Since $a \geq \frac{n}{8}$, the weighted sum on the right-hand side is at most $\left( \frac{33}{32} - c \right) n + o(n)$, for $c = \frac{1}{256}$. This completes the proof of Theorem \ref{thm:cases2to4}.
\end{proof}

The third and final step in the proof of Theorem \ref{thm:B} is to deal with
Case 1 from Section \ref{sect:kfour}. Following the derivation of
(\ref{eq:E1bound}), in order to have
$|S| \geq \left( \frac{33}{32} - o(1) \right)n$, we must have
  \begin{eqnarray}\label{eq:E1new}
    m_P + m_N \geq \frac{1}{4} (g_P + g_N) + \frac{n}{32} - o(n).
  \end{eqnarray}
  This immediately implies $P \cap M$ and $N \cap M$ can't both be empty, for
  large enough $n$, because then we'd have $m_P = m_N = 0$.

  Without loss of generality, we assume $N \cap M$ is empty and $P \cap M$
  non-empty. Thus $m_N = 0$ and in (\ref{eq:E1new}) we must have
  \begin{eqnarray}\label{eq:mP}
    m_P \geq \frac{1}{4} (g_P + g_N) + \frac{n}{32} - o(n).
  \end{eqnarray}
  But Lemma \ref{lem:M4trade} applies, and so
  \begin{eqnarray*}
    m_P \leq \frac{1}{4} g_P + \kappa(u) + o(n) \leq \frac{1}{4} g_P + \frac{n}{32} + o(n).
  \end{eqnarray*}
  This forces $g_N = o(n)$, and thus also $q_N = o(n)$, by Lemma
  \ref{lem:Q4trade}(a). It also forces $\kappa(u) = \frac{n}{32} - o(n)$, in
  other words $a = \min (P \cap M) = \frac{n}{8} + o(n)$.

  Suppose we can show that $g_P = o(n)$. Then Lemma \ref{lem:Q4trade}(a) would
  force $q_P = o(n)$ and Lemma \ref{lem:L4trade}(a) would force
  $l_P = l_N = o(n)$.
  It would just remain to show that $P \cap M$ contains at most $o(n)$ elements
  above $\frac{5n}{32}$. As in Section \ref{sect:kfour}, let
  $\theta = \frac{1}{4} \left( \frac{n}{2} + a \right) = \frac{5n}{32} + o(n)$.
  The map $y \mapsto 4y - a$ is an injection from elements of $P \cap M$ above
  $\theta$ into $P^c \cap T$. Since $g_P = o(n)$, we would be done.

  So it just remains to prove that $g_P = o(n)$. Let $\sigma_1$ and $\sigma_2$
  be the smallest and largest elements of $P \cap T$, respectively. Divide
  $P^c \cap T$ into three subsets
  \begin{eqnarray*}
    P_1 = \left\{ x : \frac{n}{2} < x < \sigma_1 \right\}, \;\; P_2 = \left\{ x \in P^c : \sigma_1 < x < \sigma_2 \right\}, \;\; P_3 = \left\{ x : \sigma_2 < x \leq n \right\},
  \end{eqnarray*}
  and let $g_i = |P_i|$, $i = 1,\,2,\,3$. It suffices to show $g_i = o(n)$ for
  each $i$.
  \\
  \\
  {\sc Step 1:} $g_3 = o(n)$.
  \\
  \par If we examine the proof of Lemma \ref{lem:M4trade}, applied now to $P$,
  then all the elements of $4A - A$ lie below
  $4\left(\frac{n}{4}\right) - a = n-a$, and hence what we actually show in
  that proof is that
  \begin{eqnarray*}
    m_P \leq \frac{g_P - \gamma}{4} + \kappa(u) + O(n^{4/5}),
  \end{eqnarray*}
  where $\gamma$ is the number of elements of $P^c \cap T$ above $n-a$.
  Comparing with (\ref{eq:mP}), we must have $\gamma = o(n)$. But
  $n-a \leq \frac{7n}{8}$, so this forces $g_3 = o(n)$. 
  \\
  \\
  {\sc Step 2:} $l_P = o(n)$.
  \\
  \par Consider $\{\sigma_2 - 4x : x \in P \cap L\}$. This set contains $l_P$
  elements. They are all positive, since $\sigma_2 > \frac{n}{2}$ and
  $\max L \leq \frac{n}{8}$. In fact, since $\sigma_2 = n - o(n)$ by Step 1,
  all but $o(n)$ of them exceed $\frac{n}{2}$. But they are all in $N^c$, by
  Lemma \ref{lem:elem}(v). Thus, $l_P \leq g_N + o(n) = o(n)$.
  \\
  \\
  {\sc Step 3:} $g_2 = o(n)$.
  \\
  \par Note that, as well as Lemma \ref{lem:M4trade}, the derivation of
  (\ref{eq:E1bound}) used Lemmas \ref{lem:L4trade}(a) and \ref{lem:Q4trade}(a),
  and all three inequalities must be equalities up to $o(n)$ in order to have
  $|S| \geq \frac{33n}{32} - o(n)$. In particular, we must have
  \begin{eqnarray*}
    l_P + l_N \geq \frac{1}{4} (g_P + g_N) - o(n).
  \end{eqnarray*}
  But we already know that $g_N = l_P = g_3 = o(n)$, so this can be written as
  \begin{eqnarray}\label{eq:g2lower}
    l_N \geq \frac{g_1 + g_2}{4} - o(n).
  \end{eqnarray}
  Now consider $\{\sigma_1 + 4x: x \in N \cap L\}$. There are $l_N$ such
  elements and all exceed $\frac{n}{2}$ since $\sigma_1$ does so. None of them
  are in $N$, by Lemma \ref{lem:elem}(iv). Since
  $\sigma_1 = \frac{n}{2} + g_1 + O(1)$,
  \begin{eqnarray*}
    \sigma_1 + 4x > n \;\; \Leftrightarrow \;\; x > \frac{n}{8} - \frac{g_1}{4} + O(1).
  \end{eqnarray*}
  It follows that $l_N \leq g_N + \frac{g_1}{4} + O(1) = \frac{g_1}{4} + o(n)$.
  Comparing with (\ref{eq:g2lower}) yields $g_2 = o(n)$.
  \\
  \\
  {\sc Step 4:} $g_1 = o(n)$ and hence $g_P = o(n)$.
  \\
  \par We have not yet used that the inequality in Lemma \ref{lem:Q4trade}(a)
  must be tight, but will do so in this last step. Thus, we must have
  \begin{eqnarray}\label{eq:qlower}
    q_P \geq \frac{g_P}{2} - o(n) \;\; \Leftrightarrow \;\; q_P \geq \frac{g_1}{2} - o(n).
    \end{eqnarray}
    We will also employ a technical general lemma about integer sumsets. All
    the variables in the following statement are assumed to be integer-valued. 

    \begin{lemma}\label{lem:sumsetgap}
      Let $\alpha \leq \beta$ and suppose $Z \subseteq [\alpha, \, \beta]$
      misses exactly $h$ elements of the closed interval. Then
      $Z+Z \supseteq [2(\alpha + h), \, 2(\beta - h)]$.
    \end{lemma}

    \begin{proof}
      We may assume that $\beta - \alpha \geq 2h$, as otherwise the last
      interval above is empty and the lemma holds trivially. Fix a $t$ in this
      interval. The correspondence $a \leftrightarrow t-a$ pairs off elements
      of the interval $[\max \{\alpha, \, t-\beta\}, \, \min \{\beta, \, t-\alpha\}] \subseteq [\alpha,\,\beta]$, with at most one integer being paired with
      itself. This interval has $\min \{\beta - \alpha, \, t-2\alpha, \, 2\beta - t\} + 1 \geq 2h+1$ elements (since $\beta - \alpha \geq 2h$). Since $Z$
      misses only $h$ elements of $[\alpha,\,\beta]$, there must therefore be
      some $a$ for which both $a$ and $t-a$ are in $Z$, in which case
      $t = a + (t-a) \in Z + Z$.
    \end{proof}

    Let $A = P \cap T$ and $B = P \cap Q$ (same notation as in Definition
    \ref{defi:exceptional}), so $|B| = q_P$.
    $A \subseteq [\sigma_1, \, \sigma_2]$ and misses $g_2$ elements of the
    interval. By Lemma \ref{lem:sumsetgap},
    $A+A \supseteq [2\sigma_1 + 2g_2, \, 2\sigma_2 - 2g_2]$. If
    $b \in B$ then $4b \not\in A+A$, by Lemma \ref{lem:elem}(ii). Thus
    $B = B_1 \cup B_2$, where
    \begin{eqnarray*}
      B_1 = B \cap \left( \frac{n}{4}, \, \frac{\sigma_1 + g_2}{2} \right), \;\;\; B_2 = B \cap \left( \frac{\sigma_2 - g_2}{2}, \, \frac{n}{2} \right].
    \end{eqnarray*}
    Now $|B_2| \leq \frac{1}{2} ((n-\sigma_2) + g_2) + O(1) = \frac{1}{2} (g_3 + g_2) + O(1) = o(n)$. Hence, by (\ref{eq:qlower}), we must have
    $|B_1| \geq \frac{g_1}{2} - o(n)$. Noting that
    $\sigma_1 = \frac{n}{2} + g_1 + O(1)$, $B_1$ is the intersection of $B$
    with an interval of length
    $\frac{1}{2} (g_1 + g_2) + O(1) = \frac{g_1}{2} + o(n)$. Hence, $B_1$
    misses at most $o(n)$ elements of the interval. Writing the open interval
    $\left( \frac{n}{4}, \, \frac{\sigma_1 + g_2}{2} \right)$ instead as a
    closed interval $[\alpha, \, \beta]$, let $h = o(n)$ be the
    number of elements in $[\alpha, \, \beta] \backslash B$.

    Applying Lemma \ref{lem:sumsetgap} again,
    $B_1 + B_1 \supseteq [2\alpha + 2h, \, 2\beta - 2h]$. Note that
    \begin{eqnarray*}
      2\alpha + 2h = \left( \frac{n}{2} + O(1) \right) + o(n) \leq \frac{n}{2} + o(n), \\ 2\beta - 2h = \left( \sigma_1 + g_2 + O(1) \right) - o(n) \geq \frac{n}{2} + g_1 - o(n).
    \end{eqnarray*}
    If $y \in P \cap M \subseteq \left[ \frac{n}{8}, \, \frac{n}{4} \right]$,
    then $4y \not\in B_1 + B_1$, by Lemma \ref{lem:elem}(ii). All of this
    implies that $P \cap M$ contains at most $o(n)$ elements in the interval
    $\left[ \frac{n}{8}, \, \frac{n}{8} + \frac{g_1}{4} \right]$. Let $Z$ be
    the set of elements in $P \cap M$ above $\frac{n}{8} + \frac{g_1}{4}$. If
    $Z \neq \phi$, we can apply the proof of Lemma \ref{lem:M4trade} to
    conclude that
    \begin{eqnarray*}
      |Z| \leq \frac{g_P}{4} + \kappa \left( \frac{n}{8} + \frac{g_1}{4} \right) + o(n)
    \end{eqnarray*}
    and hence that the same inequality holds for $m_P$, since
    $m_P \leq |Z| + o(n)$. But, since $g_N = o(n)$, this contradicts
    (\ref{eq:mP}) unless $g_1 = o(n)$.

    This completes Step 4, and with it the full proof of Theorem \ref{thm:B}.
    
\vspace{3mm}

\section{{\bf Continuous analogues}}\label{sect:cont}

Each of the three discrete settings in which $\mathcal{L}$-free sets were
introduced in Section \ref{sect:intro}, that is, subsets of
$[1,\,n]$, $[-n,\,n]$ and $\mathbb{Z}_n$, has a natural continuous analogue,
namely measurable subsets of $[0,\,1]$, $[-1,\,1]$ and $\mathbb{R}/\mathbb{Z}$
respectively. There is a considerable literature on the first and last of
these. In Section \ref{sect:intro} we noted a fundamental difference
between them, namely that, in the modular setting, the supremum of
the measures of
$\mathcal{L}$-free subsets of $\mathbb{R}/\mathbb{Z}$ equals the
limit of the maximum densities for subsets of
$\mathbb{Z}_p$, as $p \rightarrow \infty$. This ``transference''
principle does not apply in general in the integer setting.
We noted in Section \ref{sect:intro} that
$D^3 = \frac{1}{2}$, with the odd numbers in $[1,\,n]$ forming the unique
maximum $3$-sum-free set for all
$n \gg 0$. Such a set has no continuous analogue, and indeed it was proven
in \cite{PR} that the maximum possible measure of a $3$-sum-free subset of
$[0,\,1]$ is $\frac{77}{177}$, achieved by the continuous $k=3$ analogue of
(\ref{eq:Splus}).

Just as nobody seems to have looked
specifically at $\mathcal{L}$-free subsets of $[-n,\,n]$ before, so we could
not find any literature on measurable $\mathcal{L}$-free subsets of $[-1,\,1]$.
Moreover, we could not find a restricted version of a
transference principle, valid for sets of integers, explicitly
presented anywhere, even though such a principle seems to be ``well-known
folklore''. Informally, it states that
\begin{itemize}
\item
  The supremum of the measures of $\mathcal{L}$-free sets
  in the continuous setting
  cannot exceed the maximum asymptotic density in the corresponding discrete
  setting.
\item
  If the extremal sets in the discrete setting are, up to an $o(1)$ term in
  density, unions of finitely many intervals, then the extremal sets in the
  continuous setting have the same structure, up to a null set.
\end{itemize}

For the sake of completeness, we will first present a rigorous verification of
this ``integer transference principle''. 
Propositions \ref{prop:maxC} and \ref{prop:stabC} below are the relevant
results. As an immediate corollary, we will get continuous analogues of
Theorems \ref{thm:A} and \ref{thm:B}. From now until the end of the proof of
Proposition \ref{prop:stabC}, $\mathcal{L}$ always denotes a
non-invariant equation. 

\begin{lemma}\label{lem:C1}
  Let $I \subseteq \mathbb{R}$ be an interval of finite length $l$ and let
  $n \in \mathbb{N}$. Then
  \begin{eqnarray*}
    nl - 1 \leq |\mathbb{Z} \cap nI| \leq nl + 1.
  \end{eqnarray*}
\end{lemma}

\begin{proof}
This follows immediately from the fact that $nI$ is an interval of length $nl$. \end{proof}

\begin{lemma}\label{lem:C2}
  Fix any three positive integers $a,\,b,\,c$ and
  $\mathcal{L} = \mathcal{L}_{a,\,b,\,c}$. Let $I \subseteq \mathbb{R}$ be an
  interval of finite positive length and $A \subseteq I$ a measurable
  $\mathcal{L}$-free set. Then, for every $\varepsilon > 0$, there exists a set
  $W \subseteq I$ such that
  \par (i) $W$ is a union of finitely many, pairwise disjoint intervals.
  \par (ii) $W$ is $\mathcal{L}$-free.
  \par (iii) $\mu(A \Delta W) < \varepsilon$.
\end{lemma}

\begin{proof}
  Set $s := |a| + |b| + |c|$ and fix an $\varepsilon > 0$. If
  $\mu(A) < \varepsilon$, choose $W = \phi$. Otherwise, let $K \subseteq A$ be
  a compact set such that $\mu(A \backslash K) < \frac{\varepsilon}{2}$.
  $K^3 \subseteq \mathbb{R}^3$ is compact and the function
  $f: (x,\,y,\,z) \mapsto ax-by-cz$ is continuous, so $f(K^3)$ is a compact
  set. Since $K \subseteq A$, this set doesn't contain zero and hence
  $\eta = \min \{ |t| : t \in f(K^3)\}$ exists and is strictly positive.

  For any $r > 0$, set $U_r := K + (-r,\,r)$. So the $U_r$ are open sets
  containing $K$, whose intersection is $K$. Fix an
  $r < \frac{\eta}{s}$ such that
  $\mu(U_r \backslash K) < \frac{\varepsilon}{2}$ and put $U = U_r$. The choice
  of $r$ implies that $U$ is also $\mathcal{L}$-free.

  Since $U$ is open, it is a disjoint union of countably many open intervals.
  Finitely many of these must cover the compact set $K$, say intervals
  $J_1,\,J_2,\,\dots,\,J_m$. Put $W = (J_1 \cup \dots \cup J_m) \cap I$. It is
  easily checked that $W$ has the desired properties (i)-(iii).
\end{proof}

\begin{lemma}\label{lem:C3}
  Let $I \subseteq \mathbb{R}$ be a finite length interval and $W \subseteq I$
  be an $\mathcal{L}$-free set which is a union of $m$ pairwise disjoint
  intervals. For every $n \in \mathbb{N}$, the set $S_n = \mathbb{Z} \cap nW$
  is an $\mathcal{L}$-free subset of $\mathbb{Z} \cap nI$ and
  $|S_n| \geq n \cdot \mu(W) - m$.
\end{lemma}

\begin{proof}
  The lower bound on $|S_n|$ follows from Lemma \ref{lem:C1}. If
  $(x,\,y,\,z) \in S^3$ satisfied $ax=by+cz$, then so would
  $\left( \frac{x}{n}, \, \frac{y}{n}, \, \frac{z}{n} \right) \in W^3$, which
  is impossible since $W$ is $\mathcal{L}$-free. Thus $S_n$ is also
  $\mathcal{L}$-free.
\end{proof}

For a finite length interval $I \subseteq \mathbb{R}$ and positive integer
$n \in \mathbb{N}$, let
\begin{eqnarray*}
  d_I = \sup \{\mu(A) : {\hbox{$A \subseteq I$ is measurable and $\mathcal{L}$-free}} \}, \\
  D_{I,\,n} = \max \{|S| : {\hbox{$S \subseteq \mathbb{Z} \cap nI$ and $S$ is $\mathcal{L}$-free}} \}.
\end{eqnarray*}

\begin{proposition}\label{prop:maxC}
  \begin{eqnarray}\label{eq:maxC}
    d_I \leq \liminf_{n \rightarrow \infty} \frac{D_{I,\,n}}{n}.
  \end{eqnarray}
\end{proposition}

\begin{proof}
  Let $A \subseteq I$ be measurable and $\mathcal{L}$-free. Fix
  $n \in \mathbb{N}$ and $\varepsilon > 0$. Choose $W$ as in
  Lemma \ref{lem:C2}, then $S_n$ as in Lemma \ref{lem:C3}. We have
  \begin{eqnarray*}
    \frac{D_{I,\,n}}{n} \geq \frac{|S_n|}{n} \geq \mu(W) - \frac{m}{n} > \mu(A) - \varepsilon - \frac{m}{n}.
  \end{eqnarray*}
  First let $n \rightarrow \infty$, then let $\varepsilon \rightarrow 0$ and
  finally take the supremum over $A$.
  \end{proof}

\begin{proposition}\label{prop:stabC}
  Let $I \subseteq \mathbb{R}$ as before. Let $d > 0$ and let
  $B_1,\,B_2,\,\dots,\,B_r \subseteq I$ be sets, each of which is a union of
  finitely many intervals. Assume the following holds:
  \\
  \par ``If $n_j \rightarrow \infty$ and $S_j \subseteq \mathbb{Z} \cap n_j I$ are $\mathcal{L}$-free with $|S_j| \geq dn_j - o(n_j)$, then $\min_{i} |S_j \Delta (\mathbb{Z} \cap n_j B_i)| = o(n_j)$''.
  \\
  \\
  Then every measurable $\mathcal{L}$-free set $A \subseteq I$ with
  $\mu(A) \geq d$ satisfies $\mu(A \Delta B_i) = 0$ for some $i$.
\end{proposition}

\begin{proof}
  Fix a measurable $\mathcal{L}$-free $A \subseteq I$ such that
  $\mu(A) \geq d$. For each $j \geq 1$, let $\varepsilon_j = \frac{1}{j}$. Lemma
  \ref{lem:C2} yields an $\mathcal{L}$-free $W_j \subseteq I$, which is a union
  of finitely many intervals, say $m_j$ of them, and such that
  $\mu(A \Delta W_j) < \varepsilon_j$. In particular,
  $\mu(W_j) > d - \varepsilon_j$.

  Choose integers $n_j \geq \frac{m_j}{\varepsilon_j}$, so that
  $n_j \rightarrow \infty$. Let $S_j = \mathbb{Z} \cap n_j W_j$. By Lemma
  \ref{lem:C3}, $S_j$ is $\mathcal{L}$-free and
  $|S_j| \geq n_j \mu(W_j) - m_j > n_j(d-\varepsilon_j) - \varepsilon_j n_j = dn_j - o(n_j)$. 
  The assumption of the proposition thus implies that, for each $j$ there
  exists an $i = i(j) \in \{1,\,\dots,\,r\}$ such that 
  $|S_j \Delta (\mathbb{Z} \cap n_j B_i)| = o(n_j)$.
  There is at least one $i$ such that
  this holds for $B_i$ and an infinite sequence of indices $j$. Fix such an
  $i$. Then, for any corresponding $j$, 
  \begin{eqnarray*}
    S_j \Delta (\mathbb{Z} \cap n_j B_i) = (\mathbb{Z} \cap n_j W_j) \Delta (\mathbb{Z} \cap n_j B_i) = \mathbb{Z} \cap n_j (B_i \Delta W_j).
  \end{eqnarray*}
  Since $B_i$ is a union of finitely many intervals and
  $n_j \geq \frac{m_j}{\varepsilon_j}$, it easily follows from
  Lemma \ref{lem:C1} that $\mu(W_j \Delta B_i) \rightarrow 0$, as $j \rightarrow \infty$. But
  \begin{eqnarray*}
    \mu(A \Delta B_i) \leq \mu(A \Delta W_j) + \mu(W_j \Delta B_i)
  \end{eqnarray*}
  and, since both terms on the right-hand side tend to zero as
  $j \rightarrow \infty$, we conclude that $\mu(A \Delta B_i) = 0$. 
\end{proof}

\begin{remark}\label{rem:modularC}
  Note that the proof of Lemma \ref{lem:C3} utilizes the homogeneity of
  $\mathcal{L}$. It is readily checked that the whole of the above analysis
  goes through for any homogeneous, non-invariant linear equation in any number
  of integer variables. 
\end{remark}

The continuous analogues of the sets in (\ref{eq:S}) and (\ref{eq:S4}) are,
respectively,
\begin{eqnarray}\label{eq:Scont}
  S_{k} = \pm \left( \frac{2}{k}, \, 1 \right], \;\;\; \mathcal{S}_4 = \pm \left( \frac{1}{2}, \, 1 \right] \cup \left[ \frac{1}{8}, \, \frac{5}{32} \right].
\end{eqnarray}

\begin{theorem}\label{thm:D}
  Let $A$ be a measurable $k$-sum-free subset of $[-1,\,1]$.
  \par (i) If $k \geq 5$, then $\mu(A) \leq 2-\frac{4}{k}$ and, if $\mu(A) = 2 - \frac{4}{k}$, then $\mu (A \Delta S_k) = 0$.
  \par (ii) If $k=4$, then $\mu(A) \leq \frac{33}{32}$ and, if $\mu(A) = \frac{33}{32}$, then $\mu(A \Delta \mathcal{S}_4) = 0$ or $\mu(A \Delta -\mathcal{S}_4) = 0$.
\end{theorem}

\begin{proof}
  One applies Propositions \ref{prop:maxC} and \ref{prop:stabC} to
  $I = [-1,\,1]$, the latter with $d = 2 - \frac{4}{k}$, $B_1 = S_k$ for
  $k \geq 5$, and with $d = \frac{33}{32}$, $B_{1,\,2} = \pm \mathcal{S}_4$ for
  $k=4$. If one examines the proofs in previous sections, one may easily verify
  that the assumption in Proposition \ref{prop:stabC} holds. 
\end{proof}

This leaves $k\leq 3$ in the symmetric, continuous setting. For $k=1$, it's a
trivial consequence of results already in the literature (or from Proposition
\ref{prop:maxC}) that the maximum possible measure is $1$ and achieved by
$\pm \left( \frac{1}{2}, \, 1 \right]$. For $k=2$, as for any invariant
equation $\mathcal{L}$, it is likewise clear that a measurable
$\mathcal{L}$-free set must have measure zero.

For $k=3$, the known result for $[0,\,1]$ mentioned at the beginning of this
section \cite{PR} suggests the possibility of interesting
behaviour also for $[-1,\,1]$. We have not been able to determine the maximum
measure of a $3$-sum-free
subset of $[-1,\,1]$, but the evidence points to structures for the
extremal sets which differ both from the structure when $k=4$ and from that
in the case of $[0,\,1]$. Let
\begin{eqnarray}\label{eq:Bs}
  B_1 = \pm \left( \frac{2}{3}, \, 1 \right) \cup \pm \left( \frac{1}{4}, \, \frac{11}{36} \right), \\
  B_2 = \pm \left( \frac{2}{3}, \, 1 \right) \cup \pm \left( \frac{1}{9}, \, \frac{1}{6} \right).
\end{eqnarray}
Note that $\mu(B_1) = \mu(B_2) = \frac{7}{9}$, each $B_i$ is a symmetric set $(B_i = -B_i$) and $\mu(B_1 \Delta B_2) = \frac{2}{9}$, so these are two ``genuinely different'' structures, even if both are based on the same idea of full edge intervals supplemented by interior intervals whose sizes are maximized while not creating solutions to $3x=y+z$. 

\begin{conjecture}\label{conj:D}
  Let $A$ be a measurable, $3$-sum-free subset of $[-1,\,1]$. Then $\mu(A) \leq \frac{7}{9}$. Moreover, if $\mu(A) = \frac{7}{9}$, then $\mu(A \Delta B_1) = 0$ or $\mu(A \Delta B_2) = 0$.
\end{conjecture}

\vspace{3mm}

\section{{\bf AI statement}}\label{sect:AI}

The first author had the idea many years ago to study $k$-sum-free subsets of
$[-n,\,n]$ and conjectured, but was unable to prove, that the set
$\mathcal{S}_{4,\,n}$ achieved maximum density. The proofs in this paper are
essentially the work of \texttt{Claude Opus 5.5}, upon prompting by the
authors, starting with an explicit suggestion to try to prove or disprove that
conjecture. \texttt{Claude} also performed all computations, including for
the various certificates and for the numerical evidence to
support Conjectures \ref{conj:C} and \ref{conj:D}.
All the proofs have been hand-checked in full, and the paper was
entirely written by the authors. It includes some
mathematical, as well as stylistic, refinements of \texttt{Claude}'s work. In
particular, we simplified the argument for $k \geq 6$ considerably, compared to
\texttt{Claude}'s initial proposal. However, the extremely intricate analysis
for $k=4$ would have been impossible without AI assistance. Of course one may
ask whether such an intricate analysis really is necessary to prove Theorem
\ref{thm:B} or even to attack Conjectures \ref{conj:C} and \ref{conj:D}.
We leave these issues to the future.

\section*{Appendix}

\par This appendix contains the certificates (S2)-(S7), in the same format as (S1) in Section \ref{sect:kfour}.

\subsection*{Certificate (S2).} This certificate applies when $1 \leq c_P \leq 2$, $1 \leq c_N \leq 3$, $\mathcal{P}_3 = \phi$. On the right-hand sides we use the upper ends $c_P = 2$, $c_N = 3$; as the coefficient $c_O$ in (\ref{eq:DO}), (\ref{eq:DS}), (\ref{eq:DW}) we use the lower ends, i.e.~$c_N = 1$ for $X = P$ and $c_P = 1$ for $X = N$.
\begin{center}
{\scriptsize \setlength{\tabcolsep}{2pt} \renewcommand{\arraystretch}{1.4}
\resizebox{\textwidth}{!}{%
\begin{tabular}{|c|c|c|ccccccc|ccccccc|c|} \hline
Ineq. & $X$ & Weight & $l_P$ & $\mu_{1,\,P}$ & $\mu_{2,\,P}$ & $\mu_{3,\,P}$ & $q_{1,\,P}$ & $\rho_P$ & $\tau_P$ & $l_N$ & $\mu_{1,\,N}$ & $\mu_{2,\,N}$ & $\mu_{3,\,N}$ & $q_{1,\,N}$ & $\rho_N$ & $\tau_N$ & Right-hand side \\ \hline
(\ref{eq:l4a}) &  & $1$ & $1$ &  &  &  &  &  & $\frac{1}{4}$ & $1$ &  &  &  &  &  & $\frac{1}{4}$ & $\frac{n}{4}$ \\
(\ref{eq:mu4}) & $P$ & $1$ &  &  &  & $1$ &  &  &  &  &  &  &  &  &  &  & $0$ \\
(\ref{eq:Q4qr}) & $P$ & $1$ &  &  &  &  & $1$ & $1$ & $\frac{1}{2}$ &  &  &  &  &  &  &  & $\frac{n}{4}$ \\
(\ref{eq:DS}) & $P$ & $\frac{3}{16}$ &  & $5$ & $4$ &  &  &  & $1$ &  & $1$ & $1$ & $1$ &  &  &  & $\frac{3n}{4} - \frac{3a}{2}$ \\
(\ref{eq:DW}) & $P$ & $\frac{1}{16}$ &  & $1$ & $4$ &  &  &  & $1$ &  & $1$ & $1$ & $1$ &  &  &  & $\frac{5n}{8} - \frac{3a}{4}$ \\
(\ref{eq:Q4qr}) & $N$ & $\frac{9}{8}$ &  &  &  &  &  &  &  &  &  &  &  & $1$ & $1$ & $\frac{1}{2}$ & $\frac{n}{4}$ \\
(\ref{eq:WS}) & $N$ & $\frac{3}{16}$ &  &  &  &  &  &  &  &  & $4$ & $4$ & $4$ &  &  & $1$ & $\frac{7n}{8} - \frac{9b}{4}$ \\
\hline
Sum &  &  & $1$ & $1$ & $1$ & $1$ & $1$ & $1$ & $1$ & $1$ & $1$ & $1$ & $1$ & $\frac{9}{8}$ & $\frac{9}{8}$ & $1$ & $\frac{9n}{8} - \frac{21a}{64} - \frac{27b}{64}$ \\ \hline
\end{tabular}}}
\end{center}
At the five vertices of $\mathcal{R}$, listed as above, the right-hand side of the last row equals $\frac{33765n}{32768}$, $\frac{33771n}{32768}$, $\frac{521n}{512}$, $n$, $\frac{519n}{512}$. The maximum is $\frac{33771n}{32768} < \left( \frac{33}{32} - \frac{1}{2048} \right) n$.

\subsection*{Certificate (S3).} This certificate applies when $1 \leq c_P \leq 2$, $3 \leq c_N \leq 4$, $\mathcal{P}_3 \neq \phi$, $\mathcal{N}_3 = \phi$. On the right-hand sides we use the upper ends $c_P = 2$, $c_N = 4$; as the coefficient $c_O$ in (\ref{eq:DO}), (\ref{eq:DS}), (\ref{eq:DW}) we use the lower ends, i.e.~$c_N = 3$ for $X = P$ and $c_P = 1$ for $X = N$.
\begin{center}
{\scriptsize \setlength{\tabcolsep}{2pt} \renewcommand{\arraystretch}{1.4}
\resizebox{\textwidth}{!}{%
\begin{tabular}{|c|c|c|ccccccc|ccccccc|c|} \hline
Ineq. & $X$ & Weight & $l_P$ & $\mu_{1,\,P}$ & $\mu_{2,\,P}$ & $\mu_{3,\,P}$ & $q_{1,\,P}$ & $\rho_P$ & $\tau_P$ & $l_N$ & $\mu_{1,\,N}$ & $\mu_{2,\,N}$ & $\mu_{3,\,N}$ & $q_{1,\,N}$ & $\rho_N$ & $\tau_N$ & Right-hand side \\ \hline
(\ref{eq:l4a}) &  & $\frac{23}{24}$ & $1$ &  &  &  &  &  & $\frac{1}{4}$ & $1$ &  &  &  &  &  & $\frac{1}{4}$ & $\frac{n}{4}$ \\
(\ref{eq:mu2}) & $P$ & $\frac{1}{96}$ &  &  & $1$ &  &  &  &  &  &  &  &  &  &  &  & $\frac{a}{2} + \frac{b}{8} - \frac{n}{16}$ \\
(\ref{eq:rho}) & $P$ & $\frac{1}{16}$ &  &  &  &  &  & $1$ &  &  &  &  &  &  &  &  & $4b - \frac{n}{2}$ \\
(\ref{eq:Q4qr}) & $P$ & $\frac{15}{16}$ &  &  &  &  & $1$ & $1$ & $\frac{1}{2}$ &  &  &  &  &  &  &  & $\frac{n}{4}$ \\
(\ref{eq:R}) & $P$ & $\frac{1}{24}$ &  &  &  &  &  &  & $1$ & $1$ & $1$ & $1$ & $1$ & $1$ &  &  & $\frac{n}{2}$ \\
(\ref{eq:DS}) & $P$ & $\frac{1}{12}$ &  & $7$ & $4$ &  &  &  & $1$ &  & $1$ & $1$ & $1$ &  &  &  & $\frac{3n}{4} - \frac{3a}{2}$ \\
(\ref{eq:WW}) & $P$ & $\frac{3}{32}$ & $1$ & $1$ & $4$ & $4$ &  &  & $1$ &  &  &  &  &  &  &  & $\frac{5n}{8} - \frac{3a}{4}$ \\
(\ref{eq:U1}) & $P$ & $\frac{7}{96}$ & $1$ & $1$ & $1$ & $5$ &  &  & $1$ &  &  &  &  &  &  &  & $\frac{n}{2}$ \\
(\ref{eq:mu4}) & $N$ & $1$ &  &  &  &  &  &  &  &  &  &  & $1$ &  &  &  & $0$ \\
(\ref{eq:rho}) & $N$ & $\frac{1}{24}$ &  &  &  &  &  &  &  &  &  &  &  &  & $1$ &  & $4a - \frac{n}{2}$ \\
(\ref{eq:Q4qr}) & $N$ & $\frac{23}{24}$ &  &  &  &  &  &  &  &  &  &  &  & $1$ & $1$ & $\frac{1}{2}$ & $\frac{n}{4}$ \\
(\ref{eq:R}) & $N$ & $\frac{1}{16}$ & $1$ & $1$ & $1$ & $1$ & $1$ &  &  &  &  &  &  &  &  & $1$ & $\frac{n}{2}$ \\
(\ref{eq:DS}) & $N$ & $\frac{1}{6}$ &  & $1$ & $1$ & $1$ &  &  &  &  & $5$ & $4$ &  &  &  & $1$ & $n - 3b$ \\
(\ref{eq:DW}) & $N$ & $\frac{5}{96}$ &  & $1$ & $1$ & $1$ &  &  &  &  & $1$ & $4$ &  &  &  & $1$ & $\frac{5n}{8} - \frac{3b}{4}$ \\
\hline
Sum &  &  & $\frac{19}{16}$ & $\frac{33}{32}$ & $\frac{103}{96}$ & $\frac{49}{48}$ & $1$ & $1$ & $1$ & $1$ & $\frac{97}{96}$ & $1$ & $\frac{9}{8}$ & $1$ & $1$ & $1$ & $\frac{1643n}{1536} - \frac{3a}{128} - \frac{221b}{768}$ \\ \hline
\end{tabular}}}
\end{center}
At the five vertices of $\mathcal{R}$, listed as above, the right-hand side of the last row equals $\frac{405091n}{393216}$, $\frac{67549n}{65536}$, $\frac{2109n}{2048}$, $\frac{2345n}{2304}$, $\frac{9389n}{9216}$. The maximum is $\frac{67549n}{65536} < \left( \frac{33}{32} - \frac{1}{2048} \right) n$.

\subsection*{Certificate (S4).} This certificate applies when $1 \leq c_P \leq 3$, $c_N = 3$, $\mathcal{P}_3 \neq \phi$, $\mathcal{N}_3 \neq \phi$. On the right-hand sides we use the upper ends $c_P = 3$, $c_N = 3$; as the coefficient $c_O$ in (\ref{eq:DO}), (\ref{eq:DS}), (\ref{eq:DW}) we use the lower ends, i.e.~$c_N = 3$ for $X = P$ and $c_P = 1$ for $X = N$.
\begin{center}
{\scriptsize \setlength{\tabcolsep}{2pt} \renewcommand{\arraystretch}{1.4}
\resizebox{\textwidth}{!}{%
\begin{tabular}{|c|c|c|ccccccc|ccccccc|c|} \hline
Ineq. & $X$ & Weight & $l_P$ & $\mu_{1,\,P}$ & $\mu_{2,\,P}$ & $\mu_{3,\,P}$ & $q_{1,\,P}$ & $\rho_P$ & $\tau_P$ & $l_N$ & $\mu_{1,\,N}$ & $\mu_{2,\,N}$ & $\mu_{3,\,N}$ & $q_{1,\,N}$ & $\rho_N$ & $\tau_N$ & Right-hand side \\ \hline
(\ref{eq:l4a}) &  & $\frac{23}{32}$ & $1$ &  &  &  &  &  & $\frac{1}{4}$ & $1$ &  &  &  &  &  & $\frac{1}{4}$ & $\frac{n}{4}$ \\
(\ref{eq:rho}) & $P$ & $\frac{3}{32}$ &  &  &  &  &  & $1$ &  &  &  &  &  &  &  &  & $4b - \frac{n}{2}$ \\
(\ref{eq:Q4qr}) & $P$ & $\frac{59}{64}$ &  &  &  &  & $1$ & $1$ & $\frac{1}{2}$ &  &  &  &  &  &  &  & $\frac{n}{4}$ \\
(\ref{eq:R}) & $P$ & $\frac{5}{64}$ &  &  &  &  &  &  & $1$ & $1$ & $1$ & $1$ & $1$ & $1$ &  &  & $\frac{n}{2}$ \\
(\ref{eq:DS}) & $P$ & $\frac{5}{64}$ &  & $7$ & $4$ &  &  &  & $1$ &  & $1$ & $1$ & $1$ &  &  &  & $\frac{7n}{8} - \frac{9a}{4}$ \\
(\ref{eq:WS}) & $P$ & $\frac{1}{32}$ & $1$ & $4$ & $4$ & $4$ &  &  & $1$ &  &  &  &  &  &  &  & $\frac{7n}{8} - \frac{9a}{4}$ \\
(\ref{eq:WW}) & $P$ & $\frac{7}{64}$ & $1$ & $1$ & $4$ & $4$ &  &  & $1$ &  &  &  &  &  &  &  & $\frac{5n}{8} - \frac{3a}{4}$ \\
(\ref{eq:U1}) & $P$ & $\frac{1}{16}$ & $1$ & $1$ & $1$ & $5$ &  &  & $1$ &  &  &  &  &  &  &  & $\frac{n}{2}$ \\
(\ref{eq:mu1}) & $N$ & $\frac{1}{64}$ &  &  &  &  &  &  &  &  & $1$ &  &  &  &  &  & $\frac{3n}{32} - \frac{9b}{16}$ \\
(\ref{eq:rho}) & $N$ & $\frac{5}{64}$ &  &  &  &  &  &  &  &  &  &  &  &  & $1$ &  & $4a - \frac{n}{2}$ \\
(\ref{eq:Q4qr}) & $N$ & $\frac{59}{64}$ &  &  &  &  &  &  &  &  &  &  &  & $1$ & $1$ & $\frac{1}{2}$ & $\frac{n}{4}$ \\
(\ref{eq:R}) & $N$ & $\frac{5}{64}$ & $1$ & $1$ & $1$ & $1$ & $1$ &  &  &  &  &  &  &  &  & $1$ & $\frac{n}{2}$ \\
(\ref{eq:DS}) & $N$ & $\frac{5}{64}$ &  & $1$ & $1$ & $1$ &  &  &  &  & $5$ & $4$ &  &  &  & $1$ & $\frac{7n}{8} - \frac{9b}{4}$ \\
(\ref{eq:WS}) & $N$ & $\frac{5}{64}$ &  &  &  &  &  &  &  & $1$ & $4$ & $4$ & $4$ &  &  & $1$ & $\frac{7n}{8} - \frac{9b}{4}$ \\
(\ref{eq:WW}) & $N$ & $\frac{1}{16}$ &  &  &  &  &  &  &  & $1$ & $1$ & $4$ & $4$ &  &  & $1$ & $\frac{5n}{8} - \frac{3b}{4}$ \\
(\ref{eq:U1}) & $N$ & $\frac{1}{16}$ &  &  &  &  &  &  &  & $1$ & $1$ & $1$ & $5$ &  &  & $1$ & $\frac{n}{2}$ \\
\hline
Sum &  &  & $1$ & $1$ & $\frac{35}{32}$ & $\frac{33}{32}$ & $1$ & $\frac{65}{64}$ & $1$ & $1$ & $1$ & $\frac{35}{32}$ & $\frac{33}{32}$ & $1$ & $1$ & $1$ & $\frac{2123n}{2048} - \frac{a}{64} - \frac{33b}{1024}$ \\ \hline
\end{tabular}}}
\end{center}
At the five vertices of $\mathcal{R}$, listed as above, the right-hand side of the last row equals $\frac{540319n}{524288}$, $\frac{33771n}{32768}$, $\frac{25313n}{24576}$, $\frac{395n}{384}$, $\frac{527n}{512}$. The maximum is $\frac{33771n}{32768} < \left( \frac{33}{32} - \frac{1}{2048} \right) n$.

\subsection*{Certificate (S5).} This certificate applies when $1 \leq c_P \leq 4$, $c_N = 2$, $\mathcal{P}_3 \neq \phi$, $\mathcal{N}_3 \neq \phi$. On the right-hand sides we use the upper ends $c_P = 4$, $c_N = 2$; as the coefficient $c_O$ in (\ref{eq:DO}), (\ref{eq:DS}), (\ref{eq:DW}) we use the lower ends, i.e.~$c_N = 2$ for $X = P$ and $c_P = 1$ for $X = N$.
\begin{center}
{\scriptsize \setlength{\tabcolsep}{2pt} \renewcommand{\arraystretch}{1.4}
\resizebox{\textwidth}{!}{%
\begin{tabular}{|c|c|c|ccccccc|ccccccc|c|} \hline
Ineq. & $X$ & Weight & $l_P$ & $\mu_{1,\,P}$ & $\mu_{2,\,P}$ & $\mu_{3,\,P}$ & $q_{1,\,P}$ & $\rho_P$ & $\tau_P$ & $l_N$ & $\mu_{1,\,N}$ & $\mu_{2,\,N}$ & $\mu_{3,\,N}$ & $q_{1,\,N}$ & $\rho_N$ & $\tau_N$ & Right-hand side \\ \hline
(\ref{eq:l4a}) &  & $\frac{37}{48}$ & $1$ &  &  &  &  &  & $\frac{1}{4}$ & $1$ &  &  &  &  &  & $\frac{1}{4}$ & $\frac{n}{4}$ \\
(\ref{eq:mu2}) & $P$ & $\frac{1}{96}$ &  &  & $1$ &  &  &  &  &  &  &  &  &  &  &  & $a + \frac{b}{4} - \frac{n}{8}$ \\
(\ref{eq:rho}) & $P$ & $\frac{7}{96}$ &  &  &  &  &  & $1$ &  &  &  &  &  &  &  &  & $4b - \frac{n}{2}$ \\
(\ref{eq:Q4qr}) & $P$ & $\frac{89}{96}$ &  &  &  &  & $1$ & $1$ & $\frac{1}{2}$ &  &  &  &  &  &  &  & $\frac{n}{4}$ \\
(\ref{eq:R}) & $P$ & $\frac{1}{32}$ &  &  &  &  &  &  & $1$ & $1$ & $1$ & $1$ & $1$ & $1$ &  &  & $\frac{n}{2}$ \\
(\ref{eq:DS}) & $P$ & $\frac{11}{96}$ &  & $6$ & $4$ &  &  &  & $1$ &  & $1$ & $1$ & $1$ &  &  &  & $n - 3a$ \\
(\ref{eq:WW}) & $P$ & $\frac{5}{48}$ & $1$ & $1$ & $4$ & $4$ &  &  & $1$ &  &  &  &  &  &  &  & $\frac{5n}{8} - \frac{3a}{4}$ \\
(\ref{eq:U1}) & $P$ & $\frac{3}{32}$ & $1$ & $1$ & $1$ & $5$ &  &  & $1$ &  &  &  &  &  &  &  & $\frac{n}{2}$ \\
(\ref{eq:rho}) & $N$ & $\frac{1}{32}$ &  &  &  &  &  &  &  &  &  &  &  &  & $1$ &  & $4a - \frac{n}{2}$ \\
(\ref{eq:Q4qr}) & $N$ & $\frac{31}{32}$ &  &  &  &  &  &  &  &  &  &  &  & $1$ & $1$ & $\frac{1}{2}$ & $\frac{n}{4}$ \\
(\ref{eq:R}) & $N$ & $\frac{7}{96}$ & $1$ & $1$ & $1$ & $1$ & $1$ &  &  &  &  &  &  &  &  & $1$ & $\frac{n}{2}$ \\
(\ref{eq:DS}) & $N$ & $\frac{1}{24}$ &  & $1$ & $1$ & $1$ &  &  &  &  & $5$ & $4$ &  &  &  & $1$ & $\frac{3n}{4} - \frac{3b}{2}$ \\
(\ref{eq:WS}) & $N$ & $\frac{7}{48}$ &  &  &  &  &  &  &  & $1$ & $4$ & $4$ & $4$ &  &  & $1$ & $\frac{3n}{4} - \frac{3b}{2}$ \\
(\ref{eq:WW}) & $N$ & $\frac{1}{32}$ &  &  &  &  &  &  &  & $1$ & $1$ & $4$ & $4$ &  &  & $1$ & $\frac{5n}{8} - \frac{3b}{4}$ \\
(\ref{eq:U1}) & $N$ & $\frac{1}{32}$ &  &  &  &  &  &  &  & $1$ & $1$ & $1$ & $5$ &  &  & $1$ & $\frac{n}{2}$ \\
\hline
Sum &  &  & $\frac{25}{24}$ & $1$ & $\frac{35}{32}$ & $1$ & $1$ & $1$ & $1$ & $\frac{97}{96}$ & $1$ & $\frac{101}{96}$ & $\frac{97}{96}$ & $1$ & $1$ & $1$ & $\frac{205n}{192} - \frac{55a}{192} - \frac{b}{96}$ \\ \hline
\end{tabular}}}
\end{center}
At the five vertices of $\mathcal{R}$, listed as above, the right-hand side of the last row equals $\frac{16885n}{16384}$, $\frac{101257n}{98304}$, $\frac{2347n}{2304}$, $\frac{391n}{384}$, $\frac{4747n}{4608}$. The maximum is $\frac{16885n}{16384} < \left( \frac{33}{32} - \frac{1}{2048} \right) n$.

\subsection*{Certificate (S6).} This certificate applies when $1 \leq c_P \leq 4$, $3 \leq c_N \leq 4$, $\mathcal{P}_3 = \phi$. On the right-hand sides we use the upper ends $c_P = 4$, $c_N = 4$; as the coefficient $c_O$ in (\ref{eq:DO}), (\ref{eq:DS}), (\ref{eq:DW}) we use the lower ends, i.e.~$c_N = 3$ for $X = P$ and $c_P = 1$ for $X = N$.
\begin{center}
{\scriptsize \setlength{\tabcolsep}{2pt} \renewcommand{\arraystretch}{1.4}
\resizebox{\textwidth}{!}{%
\begin{tabular}{|c|c|c|ccccccc|ccccccc|c|} \hline
Ineq. & $X$ & Weight & $l_P$ & $\mu_{1,\,P}$ & $\mu_{2,\,P}$ & $\mu_{3,\,P}$ & $q_{1,\,P}$ & $\rho_P$ & $\tau_P$ & $l_N$ & $\mu_{1,\,N}$ & $\mu_{2,\,N}$ & $\mu_{3,\,N}$ & $q_{1,\,N}$ & $\rho_N$ & $\tau_N$ & Right-hand side \\ \hline
(\ref{eq:l4a}) &  & $\frac{11}{12}$ & $1$ &  &  &  &  &  & $\frac{1}{4}$ & $1$ &  &  &  &  &  & $\frac{1}{4}$ & $\frac{n}{4}$ \\
(\ref{eq:mu4}) & $P$ & $1$ &  &  &  & $1$ &  &  &  &  &  &  &  &  &  &  & $0$ \\
(\ref{eq:rho}) & $P$ & $\frac{1}{8}$ &  &  &  &  &  & $1$ &  &  &  &  &  &  &  &  & $4b - \frac{n}{2}$ \\
(\ref{eq:Q4qr}) & $P$ & $\frac{7}{8}$ &  &  &  &  & $1$ & $1$ & $\frac{1}{2}$ &  &  &  &  &  &  &  & $\frac{n}{4}$ \\
(\ref{eq:R}) & $P$ & $\frac{1}{12}$ &  &  &  &  &  &  & $1$ & $1$ & $1$ & $1$ & $1$ & $1$ &  &  & $\frac{n}{2}$ \\
(\ref{eq:DO}) & $P$ & $\frac{1}{32}$ &  & $3$ &  &  &  &  & $1$ &  & $1$ & $1$ & $1$ &  &  &  & $\frac{n}{2}$ \\
(\ref{eq:DS}) & $P$ & $\frac{1}{32}$ &  & $7$ & $4$ &  &  &  & $1$ &  & $1$ & $1$ & $1$ &  &  &  & $n - 3a$ \\
(\ref{eq:DW}) & $P$ & $\frac{3}{16}$ &  & $3$ & $4$ &  &  &  & $1$ &  & $1$ & $1$ & $1$ &  &  &  & $\frac{5n}{8} - \frac{3a}{4}$ \\
(\ref{eq:rho}) & $N$ & $\frac{5}{96}$ &  &  &  &  &  &  &  &  &  &  &  &  & $1$ &  & $4a - \frac{n}{2}$ \\
(\ref{eq:Q4qr}) & $N$ & $\frac{23}{24}$ &  &  &  &  &  &  &  &  &  &  &  & $1$ & $1$ & $\frac{1}{2}$ & $\frac{n}{4}$ \\
(\ref{eq:R}) & $N$ & $\frac{1}{8}$ & $1$ & $1$ & $1$ & $1$ & $1$ &  &  &  &  &  &  &  &  & $1$ & $\frac{n}{2}$ \\
(\ref{eq:WS}) & $N$ & $\frac{1}{6}$ &  &  &  &  &  &  &  &  & $4$ & $4$ & $4$ &  &  & $1$ & $n - 3b$ \\
\hline
Sum &  &  & $\frac{25}{24}$ & $1$ & $1$ & $\frac{9}{8}$ & $1$ & $1$ & $1$ & $1$ & $1$ & $1$ & $1$ & $\frac{25}{24}$ & $\frac{97}{96}$ & $1$ & $\frac{397n}{384} - \frac{5a}{192}$ \\ \hline
\end{tabular}}}
\end{center}
At the five vertices of $\mathcal{R}$, listed as above, the right-hand side of the last row equals $\frac{1583n}{1536}$, $\frac{33769n}{32768}$, $\frac{593n}{576}$, $\frac{593n}{576}$, $\frac{1583n}{1536}$. The maximum is $\frac{1583n}{1536} < \left( \frac{33}{32} - \frac{1}{2048} \right) n$.

\subsection*{Certificate (S7).} This certificate applies when $3 \leq c_P \leq 4$, $c_N = 4$, $\mathcal{P}_3 \neq \phi$, $\mathcal{N}_3 \neq \phi$. On the right-hand sides we use the upper ends $c_P = 4$, $c_N = 4$; as the coefficient $c_O$ in (\ref{eq:DO}), (\ref{eq:DS}), (\ref{eq:DW}) we use the lower ends, i.e.~$c_N = 4$ for $X = P$ and $c_P = 3$ for $X = N$.
\begin{center}
{\scriptsize \setlength{\tabcolsep}{2pt} \renewcommand{\arraystretch}{1.4}
\resizebox{\textwidth}{!}{%
\begin{tabular}{|c|c|c|ccccccc|ccccccc|c|} \hline
Ineq. & $X$ & Weight & $l_P$ & $\mu_{1,\,P}$ & $\mu_{2,\,P}$ & $\mu_{3,\,P}$ & $q_{1,\,P}$ & $\rho_P$ & $\tau_P$ & $l_N$ & $\mu_{1,\,N}$ & $\mu_{2,\,N}$ & $\mu_{3,\,N}$ & $q_{1,\,N}$ & $\rho_N$ & $\tau_N$ & Right-hand side \\ \hline
(\ref{eq:l4a}) &  & $\frac{47}{64}$ & $1$ &  &  &  &  &  & $\frac{1}{4}$ & $1$ &  &  &  &  &  & $\frac{1}{4}$ & $\frac{n}{4}$ \\
(\ref{eq:rho}) & $P$ & $\frac{9}{128}$ &  &  &  &  &  & $1$ &  &  &  &  &  &  &  &  & $4b - \frac{n}{2}$ \\
(\ref{eq:Q4qr}) & $P$ & $\frac{119}{128}$ &  &  &  &  & $1$ & $1$ & $\frac{1}{2}$ &  &  &  &  &  &  &  & $\frac{n}{4}$ \\
(\ref{eq:R}) & $P$ & $\frac{9}{128}$ &  &  &  &  &  &  & $1$ & $1$ & $1$ & $1$ & $1$ & $1$ &  &  & $\frac{n}{2}$ \\
(\ref{eq:DS}) & $P$ & $\frac{11}{128}$ &  & $8$ & $4$ &  &  &  & $1$ &  & $1$ & $1$ & $1$ &  &  &  & $n - 3a$ \\
(\ref{eq:WW}) & $P$ & $\frac{1}{8}$ & $1$ & $1$ & $4$ & $4$ &  &  & $1$ &  &  &  &  &  &  &  & $\frac{5n}{8} - \frac{3a}{4}$ \\
(\ref{eq:U1}) & $P$ & $\frac{9}{128}$ & $1$ & $1$ & $1$ & $5$ &  &  & $1$ &  &  &  &  &  &  &  & $\frac{n}{2}$ \\
(\ref{eq:rho}) & $N$ & $\frac{11}{128}$ &  &  &  &  &  &  &  &  &  &  &  &  & $1$ &  & $4a - \frac{n}{2}$ \\
(\ref{eq:Q4qr}) & $N$ & $\frac{119}{128}$ &  &  &  &  &  &  &  &  &  &  &  & $1$ & $1$ & $\frac{1}{2}$ & $\frac{n}{4}$ \\
(\ref{eq:R}) & $N$ & $\frac{9}{128}$ & $1$ & $1$ & $1$ & $1$ & $1$ &  &  &  &  &  &  &  &  & $1$ & $\frac{n}{2}$ \\
(\ref{eq:DS}) & $N$ & $\frac{11}{128}$ &  & $1$ & $1$ & $1$ &  &  &  &  & $7$ & $4$ &  &  &  & $1$ & $n - 3b$ \\
(\ref{eq:WS}) & $N$ & $\frac{1}{64}$ &  &  &  &  &  &  &  & $1$ & $4$ & $4$ & $4$ &  &  & $1$ & $n - 3b$ \\
(\ref{eq:WW}) & $N$ & $\frac{7}{64}$ &  &  &  &  &  &  &  & $1$ & $1$ & $4$ & $4$ &  &  & $1$ & $\frac{5n}{8} - \frac{3b}{4}$ \\
(\ref{eq:U1}) & $N$ & $\frac{9}{128}$ &  &  &  &  &  &  &  & $1$ & $1$ & $1$ & $5$ &  &  & $1$ & $\frac{n}{2}$ \\
\hline
Sum &  &  & $1$ & $\frac{133}{128}$ & $\frac{137}{128}$ & $\frac{129}{128}$ & $1$ & $1$ & $1$ & $1$ & $1$ & $\frac{137}{128}$ & $\frac{129}{128}$ & $1$ & $\frac{65}{64}$ & $1$ & $\frac{535n}{512} - \frac{a}{128} - \frac{27b}{256}$ \\ \hline
\end{tabular}}}
\end{center}
At the five vertices of $\mathcal{R}$, listed as above, the right-hand side of the last row equals $\frac{135077n}{131072}$, $\frac{67551n}{65536}$, $\frac{6331n}{6144}$, $\frac{197n}{192}$, $\frac{1051n}{1024}$. The maximum is $\frac{67551n}{65536} < \left( \frac{33}{32} - \frac{1}{2048} \right) n$.


\end{document}